\documentclass[11pt]{article}
\usepackage{amssymb,amsfonts,latexsym}
\newcommand{\cleqn}{\setcounter{equation}{0}}
\newcommand{\clth}{\setcounter{theorem}{0}}
\newcommand {\sectionnew}[1]{\section{#1}\cleqn\clth}
\newcommand{\beq}{\begin{equation}}
\newcommand{\eeq}{\end{equation}}
\newcommand{\beqa}{\begin{eqnarray}}
\newcommand{\eeqa}{\end{eqnarray}}
\newcommand{\beaa}{\begin{eqnarray*}}
\newcommand{\eaa}{\end{eqnarray*}}
\newcommand{\nn}{\hfill\nonumber}
\newcommand{\text}{\textrm}
\newcommand \nc {\newcommand}
\nc \proof {{\em{Proof.\/}} }
\nc \qed {$\Box$\hfill}
\newtheorem{theorem}{Theorem}[section]
\newtheorem{lemma}[theorem]{Lemma}
\newtheorem{definition-theorem}[theorem]{Definition-Theorem}
\newtheorem{proposition}[theorem]{Proposition}
\newtheorem{corollary}[theorem]{Corollary}
\newtheorem{definition}[theorem]{Definition}
\newtheorem{example}[theorem]{Example}
\newtheorem{remark}[theorem]{Remark}
\nc \bth[1] { \begin{theorem}\label{t#1} }
\nc \ble[1] { \begin{lemma}\label{l#1} }
\nc \bdeth[1] { \begin{definition-theorem}\label{dt#1} }
\nc \bpr[1] { \begin{proposition}\label{p#1} }
\nc \bco[1] { \begin{corollary}\label{c#1} }
\nc \bde[1] { \begin{definition}\label{d#1}\rm }
\nc \bex[1] { \begin{example}\label{e#1}\rm }
\nc \bre[1] { \begin{remark}\label{r#1}\rm }
\nc \bcon[1] { \medskip\noindent{\it{Conjecture #1}} }
\nc \bqu[1]  { \medskip\noindent{\it{Question #1}} }
\renewcommand {\eth} { \end{theorem} }
\nc {\ele} { \end{lemma} }
\nc {\edeth}{ \end{definition-theorem} }
\nc {\epr} { \end{proposition} }
\nc {\eco} { \end{corollary} }
\nc {\ede} { \end{definition} }
\nc {\eex} { \end{example} }
\nc {\ere} { \end{remark} }
\nc {\econ} {\smallskip}
\nc {\equ} {\smallskip}
\nc \eqref[1] {{\rm{(\ref{#1})}}}
\nc \thref[1]{Theorem \ref{t#1}}
\nc \leref[1]{Lemma \ref{l#1}}
\nc \prref[1]{Proposition \ref{p#1}}
\nc \coref[1]{Corollary \ref{c#1}}
\nc \deref[1]{Definition \ref{d#1}}
\nc \exref[1]{Example \ref{e#1}}
\nc \reref[1]{Remark \ref{r#1}}
\def \g  {\mathfrak{g}}
\def \a  {\mathfrak{a}}
\def \b  {\mathfrak{b}}
\def \c  {\mathfrak{z}}
\def \n  {\mathfrak{n}}
\def \m  {\mathfrak{m}}
\def \p  {\mathfrak{p}}
\def \h  {\mathfrak{h}}
\def \s  {\mathfrak{s}}
\def \l  {\mathfrak{l}}
\def \gs {\mathfrak{g}^\ast}
\def \dmf  {\mathfrak{d}}

\def \D  {\mathrm{D}}
\def \o  {\otimes}
\def \ci  {\circ}
\def \ds {\dotplus}
\def \op {\oplus}
\def \lop {\dotplus}

\def \la  {\lambda}
\def \La  {\Lambda}
\def \th {\theta}
\def \Th {\Theta}
\def \d  {\delta}
\def \sub{\subset} 
\def \sup{\supset} 
\def \C  {\mathrm{C}}
\def \pr  {\mathrm{p}}
\def \Pr  {\mathrm{p}}
\def \De  {\Delta}
\def \Sig {\Sigma}
\def \sig {\sigma}
\def \Ga {\Gamma}
\def \ga {\gamma}
\def \al {\alpha}
\def \be {\beta}
\def \Om {\Omega}

\def \i {i}

\def \e {\epsilon}
\def \O {{\mathcal O}}

\def \Rset {{\mathbb R}}
\def \Cset {{\mathbb C}}
\def \Zset {{\mathbb Z}}

\def \ra {\rightarrow}
\def \wt {\widetilde}
\def \til {\tilde}
\def \ol {\overline}
\def \wh {\widehat}

\def \hra {\hookrightarrow}
\def \st {\ast}

\def \opp { {\mathrm{op}} }
\def \Cong { {\mathrm{Cong}} }
\def \Lie { {\mathrm{Lie}} }
\def \ort { {\mathrm{ort}} }
\def \diag { {\mathrm{diag}} }
\def \tr { {\mathrm{tr}} }
\def \ord { {\mathrm{ord}} }

\def \id { {\mathrm{id}} }
\def \Ad { {\mathrm{Ad}} }

\def \rank { {\mathrm{rank}} }
\def \span { {\mathrm{span}} }

\def \Stab { {\mathrm{Stab}} }
\def \diag { {\mathrm{diag}} }

\def \Ker { {\mathrm{Ker}} }
\renewcommand \Im { {\mathrm{Im}} }
\renewcommand \max { {\mathrm{max}} }
\begin{document}
%%%%%%%%%%%%%%%%%%%%%%%%%%%%%%%%%%%%%%%%%%%%%%%%%%%%%%%%%%%%%%%%%%%%%%%
%%%%%%%%%%%%%%%%%%%%%%    Title    %%%%%%%%%%%%%%%%%%%%%%%%%%%%%%%%%%%%%%%%%%%
\title{{\LARGE\bf{
Symplectic Leaves of Complex Reductive \\ Poisson--Lie Groups}}}
\author{Milen~Yakimov 
\\{\normalsize{Department of Mathematics,}} \\
{\normalsize{University of California at Berkeley,}} \\
{\normalsize{Berkeley, CA 94720, USA}}
\thanks{E-mail address: yakimov@math.berkeley.edu}
}
\date{}
\maketitle
\begin{abstract}
All factorizable Lie bialgebra structures on complex reductive Lie
algebras were described by Belavin and Drinfeld. We classify the
symplectic leaves of the full class of corresponding connected
Poisson--Lie groups.
A formula for their dimensions is also proved.  
\end{abstract}
\setcounter{section}{-1}
%%%%%%%%%%%%%%%%%%%%   Introduction   %%%%%%%%%%%%%%%%%%%%%%%%%%%%%%%%%%%%%%%%
\sectionnew{Introduction}
%%%%%%%%%%%%%%%%%%%%%%%%%%%%%%%%%%%%%%%%%%%%%%%%%%%%%%%%%%%%%%%%%%%%%%%%%%%%%%
The structure of symplectic leaves of a Poisson manifold $M$ naturally
carries information about the geometry of the Poisson bivector field.
Its importance to dynamical systems is related to the fact that
the phase space of any Hamiltonian system on $M$ is a symplectic
leaf of $M.$ When the Poisson structure can be quantized algebraically,
it has also applications to the representation theory of the
corresponding quantized algebra of functions.
Roughly speaking the Poisson ideals of functions vanishing
on closures of symplectic leaves are deformed to primitive ideals of the
quantized algebra, but the actual relation is subtle.
In the 70's and 80's a
detailed study of the primitive ideals of universal enveloping
algebras of Lie algebras and their relation to coadjoint orbits
was done. Nice accounts to it can be found in the review
article \cite{JoIC} and the book \cite{Dix}. 
For standard (simple) Poisson--Lie groups and their quantized
algebras of functions, this relation was examined more recently
by Soibelman \cite{S2}, Hodges--Levasseur \cite{HL, HL2},
and Joseph \cite{Jo2, Jo} (see  also \cite{LS, HLT}).

In the early 80's Belavin and Drinfeld classified all nonskewsymmetric
$r$-matrices on complex simple Lie algebras. The same procedure
describes all
factorizable Lie bialgebra structures on complex reductive
Lie algebras, see \cite{Ho}. The set of symplectic
leaves on the corresponding Poisson--Lie groups is well
understood for the standard structure
\cite{DKP, HL, HKKR}, on the corresponding compact real forms 
\cite{S2, LW, S}, and in the cases of the simplest twistings
of the standard $r$-matrix by elements of the second tensor 
power of the (fixed) Cartan subalgebra
\cite{LS, HLT}. In the same time very little is known for any other
Poisson structure in this list.

The goal of this article is to present a unified description
of the sets of symplectic leaves of the complex reductive
Poisson--Lie groups associated to all $r$-matrices from
Belavin--Drinfeld's list. We also prove a formula for the
dimension of each leaf.

In the rest of this introduction we will give a rough description
of our results. Any Poisson--Lie group $G$ of the discussed class 
has two canonical Poisson--Lie subgroups $G_\pm.$ Their Lie algebras
are the spans of the second and the first components
of the associated $r$-matrix. $G_\pm$ lie inside two parabolic
subgroups $P_\pm$ of $G$ and differ from them only on a part
of a maximal torus $H \sub P_\pm.$ Let $W,$ $W_1,$ and $W_2$ denote the 
Weyl groups
of $G$ and the Levi factors  $L_1,$ $L_2$ of $P_\pm.$

The main ingredient in the classification theorem for the leaves
of $G_-$ is a choice of minimal length representative $v$ in $W$ for a coset
  {} from $W/W_1.$ (By $\dot{v}$ we will denote a representative of it
in the normalizer $N(H)$ of $H$ in $G.)$
For each such element $v$ the maximal subalgebra
of $\Lie(L_1)$ which is stable under $\Ad_{\dot{v}}$ is a reductive
subalgebra of $\g$ generated by the Cartan subalgebra $\Lie(H)$ and root
spaces
of simple roots. The corresponding reductive subgroup of
$G$ will be denoted by $G^v.$ The second main ingredient
in the classification is a choice (for each fixed $v)$
of an orbit of the twisted conjugation action of the
derived subgroup $(G^v)'$ of $G^v$ on itself defined by
\[
TC^{\dot{v}}_g (f) = g^{-1} f \Ad_{\dot{v}}(g), \; f, g \in (G^v)'   .
\]
In Sect.~4.1 we prove that this is the main part of the data which
classifies double cosets of the
dual Poisson--Lie group $G_+$ of $G_-$ in their classical double.
Sections 4.2 and 4.3 realize the passage to
symplectic leaves (dressing orbits) on $G_-.$ The set
of symplectic
leaves of $G_+$ are described in a similar way, but the corresponding
results will not be stated.

Sect.~5 deals with the general case of symplectic leaves on
the full group $G.$ The classification data is partly doubling
of the one for $G_-.$ One starts with a choice of 
a pair $(v_1, v_2)$ of minimal length representatives in $W$ for
cosets from
$W/W_1$ and $W/W_2.$ The second part of the data is a twisted 
conjugation orbit on the maximal semisimple subgroup of $L_1$ whose
Lie algebra is stable under $\Ad_{\dot{v}_1} \th^{-1} \Ad_{\dot{v}_2} 
\th.$ The map $\th: \Lie(L'_1) \ra \Lie(L'_2)$ is an isomorphism
known as
Cayley transform of the $r$-matrix $r.$ 

The first three Sections contain both a review of known facts that
we need and some new results, see in particular
Sections~1.3, 2.2, 2.3, and 3.3. 

For simplification purposes we have chosen to state and prove fully
the results for the Poisson--Lie subgroup $G_-$ of $G.$
The ones for $G$ are stated, but their proofs are mostly sketched.
Although they imply the ones for $G_-,$ their proofs compared
to the $G_-$ case contain additional complications only.

Our proofs are based on a certain inductive procedure which 
resembles (and in some sense should be considered as 
Poisson--Lie analog of) 
Duflo's construction for computing the primitive spectrum
of universal enveloping algebras of Lie algebras
which are not necessarily solvable or semisimple \cite{Du}.
The reason for this relation comes from  the fact that 
linearization of the Poisson structure on $G$ at the origin
and quantization produces the universal enveloping algebra $U(\g^\st)$ 
of the dual Lie bialgebra $\g^\st$ of $\g = \Lie(G).$
The Lie algebra $\g^\st$ is roughly speaking isomorphic
to a sum of two parabolic subalgebras of $\g$ with identified
Levi factors. It is solvable only in the cases of the simplest
twistings of the standard structure.

An explicit quantization of the Belavin--Drinfeld $r$-matrices
was obtained recently by Etingof, Schedler, and Schiffmann
\cite{EShS}. An appropriate extension of Duflo's construction
will be essential for the classification the primitive
ideals of the corresponding quantized algebras of functions.

In a subsequent publication \cite{Y} we will extend our
methods to classify the symplectic leaves of the Poisson 
homogeneous spaces for the considered groups. 

A word on terminology and notation: All double cosets will be considered
as orbits of actions of Lie groups and as such will be equipped
with the induced quotient topology. All results will refer
to this topology. Often there will be formulas that mix
direct sums of linear spaces and direct sums of Lie algebras.  
To distinguish them, the former will be denoted by $\ds$ and the
latter by $\op.$

{\flushleft{\bf{Acknowledgements}}} I am grateful to Nicolai Reshetikhin
for teaching me Poisson--Lie groups and for many stimulating discussions.
I also benefited a lot from fruitful conversations with G.~Ames,
E.~Frenkel, A.~Givental, T.~Hodges, T.~Levasseur, A.~Okounkov, 
A.~Weinstein, and J.~Wolf. 
This work was supported from the NSF grants DMS96-03239 and
DMS94-00097.   
\medskip\noindent
%%%%%%%%%%%%%%%%%%%%%%%%%%%%%%%%%%%%%%%%%%%%%%%%%%%%%%%%%%%%%%%%%%%%%%%%%%%%%%%
\sectionnew{Lie bialgebras and Poisson--Lie groups}
%%%%%%%%%%%%%%%%%%%%%%%%%%%%%%%%%%%%%%%%%%%%%%%%%%%%%%%%%%%%%%%%%%%%%%%%%%%%%%%
All Lie algebras and groups will be assumed finite dimensional.
%%%%%%%%%%%%%
\subsection{Lie bialgebras}
%%%%%%%%%%%%%
\bde{lbi} A Lie bialgebra is a  Lie algebra $(\g, \; [., .])$ equipped
with a map $\d: \g \to \wedge^2 \g,$ called Lie cobracket, such that

1. The dual map to $\d$ defines a Lie algebra structure on $\g^\st$ and

2. $\d([a, b]) = [\d(a), 1 \o b + b \o 1] 
                 + [1 \o a + a \o 1, \d(b)], \quad \forall a, b \in \g.$
\ede
The notion of Lie bialgebra is self dual. More precisely,
the dual space $\gs$ is naturally equipped with a structure of Lie
bialgebra with Lie bracket induced by $\d$ and Lie cobracket induced by
$[.,.].$ The same Lie algebra equipped with the opposite cobracket
will be denoted by $\g^{\st \opp}.$

\bde{double} The space $\g \ds \gs$ admits a canonical structure of Lie
bialgebra called classical (Drinfeld) double of $\g$ and denoted by
$\D(\g).$ It is uniquely determined from the requirements:

1. The natural embeddings $i: \: \g \hra \D(\g) = \g \ds \gs$ and 
$i^\st: \g^{\st \opp} \hra \D(\g) = \g \ds \gs$ are morphisms of
Lie bialgebras.  

2. The inner product on $\g \ds \gs$ which is equal to the 
pairing between $\g$ and $\gs$ and is trivial otherwise,
is invariant with respect to the adjoint action of $\D(\g).$
\ede

The following is an important construction for studying the structure of a
classical double. 
\bde{manin} A Manin triple is a triple of Lie algebras
$(\g, \a_+, \a_-)$ for which $\g$ is equipped with an invariant
nondegenerate inner product $(., .)$ and

1. $\a_+$ and $\a_-$ are isotropic Lie subalgebras of $\g,$ 

2. $\g = \a_+\ds \a_-.$
\ede

The double $\D(\g)$ of any Lie bialgebra $\g$ gives rise to a Manin triple 
$(\D(\g), \g, \gs)$ with inner product on $\g \ds \gs$ as defined in 
\deref{double}. Conversely, assume that $(\g, \a_+, \a_-)$ is a Manin
triple. Then $\a_+,$ $\a_-,$ and $\g$ have natural structures of Lie
bialgebras with costructures induced from the bilinear form on $\g$ by
identifying $\a_+^\st \cong \a_-,$ $\a_-^{\st \opp} \cong \a_+,$ 
and setting $\g = \a_+ \op \a_-$ as Lie coalgebras. One finds that
as bialgebras $\a_- \cong \a_+^{\st \opp}$ and $\g \cong \D(\a_+).$ 

\bde{cobl} A Lie bialgebra $(\g, \; [., .], \; \d)$ is called
coboundary if there exists an element $r \in \g \o \g$ such that
\beq
\d(a) = [ r, a \o 1 + 1 \o a], \quad \forall a \in \g.
\label{1.1}
\eeq
\ede

The map $\d$ from \eqref{1.1} defines a Lie bialgebra structure
on $\g$ if and only if the following conditions hold for $r:$
\beqa
&& r + r_{21} \in \left(\g^{\o 2}\right)^\g,
\label{1.2} \\
&& [r_{12}, r_{13}] + [r_{12}, r_{23}] + [r_{13}, r_{23}] \in 
\left(\g^{\o 3} \right)^\g.
\label{1.3} 
\eeqa
For the standard notation $r_{12},$ $r_{13},$ etc., we refer
to any book on quantum groups, e.g. \cite{CP, ESh, LS}. 
In both formulas $(.)^\g$ means $\g$-invariant part with respect to the
natural adjoint action of $\g$ on tensor powers of $\g.$ 

\bde{1.qtr} A coboundary Lie bialgebra is called quasitriangular if its
$r$-matrix satisfies the Yang--Baxter equation
\beq
[r_{12}, r_{13}] + [r_{12}, r_{23}] + [r_{13}, r_{23}] =0.
\label{1.YB}
\eeq
\ede
The classical double $\D(\g)$ of any Lie bialgebra $\g$ is a
quasitriangular bialgebra with $r$-matrix $r \in \g \o \g^\st$
representing the graph of the
trivial map $\id: \g \ra \g.$ 

For any quasitriangular Lie bialgebra $\g$ one defines the maps $r_+$
and $r_-$ $:\gs \to \g$ by
\beqa
&&r_+(x)= (x \o \id)r, \quad \forall x \in \gs,
\label{1.4} \\
&&r_-(x)= -(\id \o x)r, \quad \forall x \in \gs.
\label{1.5} 
\eeqa

\ble{1.3}The maps $r_+, \; r_-:\gs \to \g$ are morphisms of Lie
bialgebras.
\ele
\noindent
See \cite{RS, R, ESh} for a proof of this Lemma. 
Denote the images of $r_+$ and $r_-$ by $\g_+$ and $\g_-$ respectively.

\bde{1.fact} A quasitriangular Lie bialgebra $\g$ is called
factorizable \cite{RS} if
\[
r + r_{21} \in S^2 \g
\]
defines a nondegenerate inner product on $\gs.$ 
(Here $S^k V$ denotes the $k$-th symmetric power of
the vector space $V.$)
\ede

For any quasitriangular Lie bialgebra $\g,$ the spaces $\g_+$ and $\g_-$
are Lie sub-bialgebras of $\g$ and $\g_- \cong \g_+^{\st \opp}.$ Thus 
$\D(\g_+)$ as a vector space is isomorphic to $\g_+ \ds \g_-$ 
and one can define a linear map
\beq
\pi : \D(\g_+) \ra \g
\label{1.pi}
\eeq
which restricted to $\g_+$ and $\g_-$ coincides with the
natural embeddings of these algebras in $\g.$

\ble{1.pi} The map $\pi$ defined in \eqref{1.pi} is a morphism 
of Lie bialgebras.
\ele
Proofs of this Lemma can be found in \cite{RS, R, ESh}.
%%%%%%%%%%%%%%%%%%%%%%%%%%%%%%%%%%%%%%%%%
\subsection{Symplectic leaves of Poisson--Lie groups}
Recall that a symplectic leaf of a Poisson manifold
$(M, \pi)$ is a maximal connected symplectic submanifold.
Any symplectic leaf $S$ can also be described as a union of all
piecewise smooth paths starting at a given point of $S$ with segments
that are integral curves of Hamiltonian vector fields
on $M.$ In particular this implies that $M$ is a disjoint union of its
symplectic leaves. 

A Poisson--Lie group $G$ is a Lie group equipped with a
Poisson bivector field for which the multiplication map
$G \times G \ra G$
is a Poisson map. Its tangent Lie algebra $\g$ carries a natural
Lie bialgebra structure and conversely for any Lie bialgebra
$\g$ there exists a unique connected and simply connected
Poisson--Lie group whose tangent Lie bialgebra is $\g$
(see \cite{CP}). Any Poisson--Lie group $G$ with coboundary
tangent Lie bialgebra $\g$  can be equipped with the
Sklyanin Poisson--Lie structure
\beq
\pi_x = L_x (r) - R_x (r), \: x \in G
\label{sklyanin}
\eeq
where $r\in \g \o \g \equiv T_e G \o T_e G$ is
an $r$-matrix for the Lie bialgebra structure on $\g$ 
and $L_x,$ $R_x$ denote the differentials of the
right / left translations by $x \in G.$
The bivector field $\pi$ is skewsymmetric because of the 
$\Ad$-invariance of $r + r_{21},$ see \eqref{1.2}.

A double $\D(G)$ of a Poisson--Lie group $G$
is a connected Lie group with $\Lie(\D(G)) = \D(\g)$
for which the embedding $i: \g \hra \D(\g)$ can be integrated
to a (possibly not proper) embedding of Lie groups:
\[
i: G \ra \D(G).
\]
Since $\D(\g)$ is a quasitriangular Lie bialgebra, any such group 
$\D(G)$ can be equipped with the Sklyanin Poisson structure
and then $i$ becomes an embedding of Poisson--Lie groups.
By $G^r$ we will denote the (generally not closed) subgroup of $\D(G)$
with $\Lie(G^r) = i^\st(\g^\st) \sub \D(\g).$ $G^r$ will be equipped
with the topology having as a basis of open sets the connected 
components of intersections of open sets of $\D(G)$ with $G^r.$
Since $\i^\st(\g^{\st \opp})$ is a Lie sub-bialgebra of
$\D(\g),$ $G^r$ is a (nonproperly embedded) Poisson--Lie subgroup 
of $\D(G)$ with tangent Lie
bialgebra  $i^\st(\g^{\st \opp}) \cong \g^{\st \opp}.$ 
It naturally plays the role of (opposite) dual of $G.$  

The intersection
\[
\Sig \; := \; i(G) \cap G^r
\]
is a discrete subgroup of $\D(G)$ because 
$i(\g) \cap \i^\st(\g^\st) = 0.$
The right action of $\Sig$ on $i(G)$ preserves
the Poisson bivector field. By abuse of notation the projections 
of the symplectic leaves of $G$ on $i(G)/\Sig$ will
be called symplectic leaves of $i(G)/\Sig.$ 
In the case when $i(G)$ and $G^r$ are closed subgroups
of $\D(G),$ the latter are indeed the symplectic leaves of $i(G)/\Sig.$
In this case $i(G) / \Sig$ is also diffeomorphic as a Poisson manifold 
to the dense, open subset
$i(G)G^r/G^r$ of the Poisson homogeneous space $\D(G)/G^r$ of $\D(G).$ 
The projection map $G \ra i(G) / \Sig$ will be denoted by $p.$

In this setting an algebraic way to describe the symplectic leaves
of $G$ was found by Semenov-Tian-Shansky \cite{STS2} (see also
\cite{LW}). In this formulation the statement was taken from
\cite{HL}.

\bth{1.STS} 1. The symplectic leaves of $i(G) / \Sig$ are the
connected components of the intersections of double
cosets  $G^r x G^r$
in $\D(G)$ with $i(G) G^r$ projected into
$i(G) G^r/ G^r \cong i(G) / \Sig.$

2. The symplectic leaves of $G$ are connected components
of the pull backs of symplectic leaves of $i(G)/\Sig$
under the projection $p.$
\eth
In the case when $i(G)$ and $G^r$ are closed subgroups
of $\D(G),$ the left (local) action of $G^r$ on
$i(G) G^r / G^r$ can be lifted to a (local) Poisson action of $G^r$ on
$G$ using the 
the covering map $p.$ In the general case such an action is
obtained by integrating the so called dressing vector fields, 
see \cite{LW}. This action is called left dressing action of 
the (opposite) dual Poisson--Lie group $G^r$ of $G,$ on $G.$
The traditional way of stating the result of \thref{1.STS} is
that the symplectic leaves of $G$ are the orbits of the dressing action
of $G^r.$
%%%%%%%%%%%%%%%%%%%%%%%%%%%%%%%%%%
\subsection{The passage from double cosets to symplectic leaves}
%%%%%%%%
According to \thref{1.STS} the main question in passing from
$G^r$ double cosets in $\D(G)$ to symplectic leaves of
$i(G)/\Sig$ is to understand whether (and how) each such double coset
intersects the dense, open subset $i(G) G^r$ of $\D(G).$
In this Subsection we provide an answer to this question for an
arbitrary Poisson--Lie group.

\bth{int} For any connected complex Poisson--Lie group $G$
the intersection of each double coset $G^r x G^r$ in $\D(G)$
with $i(G) G^r$ is a dense, open, and connected subset
of the coset $G^r x G^r.$
\eth
Taking into account $i(\g) \ds i^\st(\g^\st) = \D(\g),$
\thref{int} can be obtained as a consequence from the following
Lemma.

\ble{inters}Let $C$ be a connected complex Lie group and 
$A,$ $B$ be two Lie subgroups of $C$ satisfying
\beq
\Lie(A) \ds \Lie(B) = \Lie(C).
\label{assume}
\eeq
Then the intersection $A x A \cap B A$ is a dense, open, and connected
subset of $A x A$ for any $x \in C.$  
\ele
\proof Eq. \eqref{assume} implies that the products $AB$ and $BA$ are
dense, open subsets of $C.$ Recall that if $U$ and $V$ are two dense, open
subsets of $C,$ then $UV =C$ (see e.g. Lemma~7.4 in \cite{Hum}). Applying
this for $U=AB$ and $V=BA$ shows 
\[
ABA =C.
\]
Thus the intersection $A x A \cap BA$ is 
a nonempty open subset of $A x A,$ for all $x \in C.$

$BA$ can be considered as the orbit
of $e \in C$ under the obvious action of $A \times B.$ 
Consider the general situation of an action of a Lie group
$F$ on a complex manifold $X$ for which there is a dense, open 
orbit $\O$ of dimension $\dim F.$ We will construct a 
global section of the anticanonical bundle $K^\st_X$ on $X$ whose 
set of 
zeros is the complement to $\O.$ Fix a linear basis
$\{f_k\}_{k=1}^N$ of $\Lie(F),$ $N =\dim F.$ The infinitesimal action
of $\Lie(F)$ gives rise to a set of vector fields $\{v_k\}_{k=1}^N$ on
$X.$
The assumption on the action of $F$ implies that the set of
vectors $\{v_1(x), \ldots, v_N(x) \} \sub T_x X$ is linearly independent
for $x \in \O$ and linearly dependent otherwise. Thus the 
global section $v_1 \wedge \ldots \wedge v_N$ of the anticanonical
bundle on $X$ vanishes exactly on the complement to the orbit $\O$
in $X.$ Applying this argument to the above mentioned action
of $A \times B$ on $C$ gives that $A x A \cap BA$ is a nonempty
open subset of $BA$ whose complement in $BA$ is the zero set
of a global section of $K^\st_C.$ This implies that it
is a dense, open subset of $BA$ which in addition is connected since
the ground field is $\Cset.$
\hfill \qed

Note also that \leref{inters} without the connectedness part
is still valid over $\Rset.$ 
%%%%%%%%%%%%%%%%%%%%%%%%%%%%%%%%%%%%%%%%%%%%%%%%%%%%%%%%%%%%%%%%%%%%%%%%%
\sectionnew{Lie algebra structure of certain doubles and duals}

This Section is devoted to the explicit description of the Lie algebra 
structure of certain doubles and duals associated to a factorizable Lie
bialgebra $\g.$ (In all statements the costructure of the involved 
bialgebras will not be considered.)    
First we recall a result of Reshetikhin and Semenov-Tian-Shansky
\cite{RS} on Lie algebra structure of $\D(\g)$
and the embeddings of $\g$ and $\gs$ in it. From it we deduce a similar
result about the Lie structure of $\D(\g_+)$
and the embeddings of $\g_+$ and $\g_-.$ 
A second derivation is also presented. Compared to the first
one, it has the the advantage that it does not rely on an initial guess of
what the structure is.  

The nondegenerate bilinear form on $\g$ induced by
$r + r_{21} \in S^2 \g$ will be used often and will be denoted
by $(.,.)_r.$
%%%%%%%%%%%%%%%%%
\subsection{Description of $\D(\g)$ and of the embeddings
$\g, \; \gs \sub \D(\g)$}
%%%%%%%%%%%%%%%%
The homomorphisms $r_\pm$ (see eqs. \eqref{1.4}--\eqref{1.5}) can be
combined to define an embedding of Lie algebras
\beq
r: \gs \ra \g \op \g
\label{1.r}
\eeq
by
\beq
r(x)= ( r_+(x), r_-(x)), \; x \in \g^\st.
\label{1.r_def}
\eeq
This map is injective due to the nondegeneracy of $r + r_{21}.$
Its image will be denoted by $\g^r \cong \gs.$

The Lie algebra $\g$ by itself can be embedded diagonally in $\g \op \g:$
\beq
d(g) = (g, g) \in \g \op \g, \; \mbox{for} \; g \in \g,
\label{1.d}
\eeq
so $d(\g) := \Im d \cong \g.$

\bpr{RS}{\em{(}}Reshetikhin--Semenov-Tian-Shansky{\em{\/)}} 
The double $\D(\g)$ of
a factorizable Lie bialgebra $\g$ is isomorphic as a Lie algebra to
$\g \op \g.$ In addition $\g^\st$ and $\g$ are embedded in it by the 
maps $r$ and $d$ defined in \eqref{1.r_def} and \eqref{1.d}.
\epr
\proof Equip $\g \op \g$ with the bilinear form 
\beq
<(x_1, x_2), (y_1, y_2)>= (x_1, y_1)_r - (x_2, y_2)_r.
\label{ggform}
\eeq
One checks that $d(\g)$ and $\g^r$ are
isotropic subspaces of $\g.$ Since $d(\g) \cap \g^r = 0$ and \\
$\dim d(\g) + \dim \g^r = 2 \dim \g$ the triple 
$(\g \op \g, d(\g), \g^r)$ is a Manin triple, which implies the statement.
\hfill \qed 

The algebra $\g^r$ can be described more explicitly in terms 
of the so called Cayley transform of $r.$ To formulate this result,
we need to introduce the algebras
\beqa
&& \m_+ = r_+ (\Ker \, r_-), \nn \\
&& \m_- = r_- (\Ker \, r_+). \nn
\eeqa
Clearly $\m_\pm$ are Lie ideals of $\g_\pm.$ The quotients
$\g_+/ \m_+$ and $\g_-/\m_-$ are isomorphic as Lie algebras
and the isomorphism is given by
\beq
\th( r_+ (x) + \m_+)= r_-(x) + \m_-, \;
\forall x \in \g^\st.
\label{defCayley}
\eeq
It is correctly defined, since $r_+ (x) \in \m_+$ implies 
$ x \in \Ker \, r_+ + \Ker \, r_-$ and thus $r_- (x) \in \m_-.$
It is also straightforward to check that $\th$ is a Lie algebra
homomorphism. The fact that it is an isomorphism is proved directly
by constructing the inverse map
$\th^{-1}: \; \g_-/ \m_- \ra \g_+/ \m_+$ by
\[
\th( r_- (x) + \m_-)= r_+(x) + \m_+, \;
\forall x \in \g^\st.
\]

The map $\th$ is the classical Cayley transform ``$\frac{f}{f-1}$ ''
of the linear map 
$f: \: \g \ra \g,$ defined by
\beq  
f = - r_- \circ j : \g \ra \g
\label{1.f}
\eeq
where $j: \: \g \ra \g^\st$ is the linear isomorphism
associated with the nondegenerate form $(.,.)_r.$
In terms of it: $r_+ \circ j = 1-f.$
The algebras $\g_\pm,$ $\m_\pm,$ are $\g_+ = \Im (f-1),$ $\g_- = \Im f,$
$\m_+ = \Ker f,$ $\m_- = \Ker (f-1)$ as needed for the definition
of Cayley transform. 

Consider also the projections:
\beq
\pr_\pm: \: \g_\pm \ra \g_\pm / \m_\pm.
\label{pr+-}
\eeq

The following result of Semenov-Tian-Shansky \cite{STS1}
describes $\g^r$ in terms of $\th.$
\bpr{1.g^r} The subalgebra $\g^r$ of $\g_+ \op \g_- \sub \g \op \g$
consists of those
$(x_+, x_-) \in \g_+ \op \g_-$ for which
\[
\th \circ \pr_+(x_+) = \pr_-(x_-).
\]
\epr
This proposition is obtained directly from the definition of $\th.$

Note that $\m_\pm$ are isotropic subalgebras of $\g_\pm$ with respect
to the bilinear form $(.,.)_r.$ Indeed, for all 
$x \in \Ker \, r_- \sub \g^\st$ and $y \in \g^\st$
we have
\beqa
&&(r_+(x), r_+(y))_r = x ( r_+(y))
\nn \\
&& = (y \op x) r= - y( r_-(x)) = 0.
\nn
\eeqa
Therefore one can equip $\g_+/ \m_+$ and $\g_-/ \m_-$ with
inner products induced by $(.,.)_r.$ In this setting the
map $\th$ is an isometry.
%%%%%%%%%%%%%%%%%%%%%%%%%%%%%%%%%%%%%%%%%%%%%
\subsection{Description of $\D(\g_+)$ and of the embeddings
$\g_+, \: \g_- \ra \D(\g_+)$}
Here we will describe the Lie structure of $\D(\g_+)$ and how $\g_+$ and
$\g_-$ are embedded in it.

\bpr{1.embed} As a Lie algebra 
\beq
\D(\g_+) \cong \g \op (\g_+/\m_+).
\label{1.doub}
\eeq
The embeddings
\[
i_\pm: \g_\pm \: \ra \D(\g_+) \cong \g \op (\g_+/\m_+)
\]
are given by the formulas
\beqa
&i_+(x)=(x, \pr_+(x)), \quad &\mbox{for} \; x \: \in \g_+,
\label{1.incl+} \\
&i_-(y)=(y, \th^{-1} \circ \pr_-(y)), \quad &\mbox{for}\; y\: \in \g_-.
\label{1.incl-}
\eeqa
\epr
\proof Equip the Lie algebra $\g \op (\g_+/\m_+)$ with the invariant
bilinear form
\beq
<(x_1, x_2), (y_1, y_2)>= (x_1, y_1)_r - (x_2, y_2)_r, \;
x_1, y_1 \in \g, \; x_2, y_2 \in \g_+/ \m_+.
\label{g+form}
\eeq
As was noticed, this is possible since $\m_+$ is an isotropic subspace of
$\g_+$ with respect to the bilinear form $(.,.)_r.$ 
Eqs. \eqref{1.incl+}--\eqref{1.incl-} define embeddings
of $\g_\pm$ in $\g \op (\g_+/\m_+).$ Their images $i_\pm(\g_\pm)$ 
are isotropic subspaces with trivial intersection.
This follows from 
\beqa
&& i_+(\g_+) = (\id \times \pr_+) d(\g_+), \nn \\
&& i_-(\g_-) = (\pr_+ \times \id) (\g^r \cap (\g_+ \op \g) ),
\nn 
\eeqa
and the fact that $d(\g)$ and $\g^r$ are isotropic subspaces
of $\g \op \g$ with respect to the inner product
\eqref{ggform}. In the second equation the elements
of $\g \op (\g_+/ \m_+)$ in the right hand side are written
as pairs in the opposite order.
Taking into account $\dim \g_+ + \dim \m_+ = \dim \g$ one gets that
$(\g \op (\g_+/ \m_+), i_+(\g_+), i_-(\g_-))$ is a Manin
triple which implies the statement. \\
\hfill \qed
%%%%%%%%%%%%%%%%%%%%%%%%%%%%%
\subsection{Another approach to the structure of $\D(\g_+)$} 
Here we will present a second approach to the results from Sect.~2.2 on
the Lie structure of $\D(\g_+).$ All results will be formulated
without proofs.

The starting point is the Lie bialgebra morphism $\pi: \D(\g_+) \ra \g,$
see \eqref{1.pi}. Its kernel 
\beq
\C \g :=
\{ (x, -x) | \; x \in \g_+ \cap \g_-\} \sub \g_+ \ds \g_-
\label{1.kerpi}
\eeq
is a Lie bialgebra ideal of $\D(\g).$
Define a map
\beq
I: \g \ra \D(\g_+)
\label{1.I}
\eeq
by
\beq
I(x) = (r_+ \circ j (x), - r_- \circ j(x))
\in \g_+ \ds \g_- \cong \D(\g_+)
\label{1.rde}
\eeq
where as before $j: \g \ra \gs$ is the identification map coming from the
nondegenerate bilinear form $(.,.)_r$ on $\g.$ 

\ble{1.rth} $I$ is an injective homomorphism of Lie algebras,
such that $r \circ I = \id_\g.$
\ele
This proposition allows us to split $\D(\g_+)$ as a direct sum 
of Lie algebras
\beq
\D(\g_+) = \Im \: I \op \C \g \cong \g \op \C \g. 
\label{dg+}
\eeq
To investigate the Lie structure of $\C \g,$ we define the maps 
\[
F_\pm: \g_\pm / \m_\pm \ra \C \g \sub \D(\g_+),
\]
by
\beqa
&&F_+(x) = (f x, -f x) \in \C \g,      \quad \forall x \: \in \Im(f-1), 
\nn \\
&&F_-(x) = ((f-1)x, -(f-1)x)\in \C \g, \quad \forall x \in \: \Im f. 
\nn
\eeqa
They are clearly well defined on the factors.
\ble{1.C} The maps $F_+$ and $F_-$ are Lie algebra isomorphisms and the
composition
\[
(F_-)^{-1} \circ F_+: \g_+/ \m_+ \ra \g_-/ \m_-
\]
is equal to the Cayley transform $\th:\g_+/ \m_+ \ra \g_-/\m_-$  
of $f.$
\ele

Combining \eqref{dg+} with this Proposition implies the isomorphism
\eqref{1.doub}. The statement about the inclusions $i_\pm$ of
$\g_\pm$ in $\D(\g_+)$ from \prref{1.embed} is obtained by computing
the composition of the embeddings $\g_\pm \hra \D(\g_+)$ with this
isomorphism. 

Finally note that when $\D(\g_+)$ is identified with
$\g \op (\g_+/ \m_+),$
the projection \\ $\pi: \: \D(\g_+) \ra \g$ 
(see \eqref{1.pi})
is just the projection on the first component.

\sectionnew{Factorizable reductive Lie bialgebras and
Poisson--Lie \\ groups}
\subsection{Belavin-Drinfeld classification}
Consider a complex reductive Lie algebra $\g.$ 
Let $\h \sub \g$ be a Cartan subalgebra of $\g.$ Choose
a set of positive roots $\De_+$ and denote by $\b_\pm$ the corresponding
positive and negative Borel subalgebras of $\g.$ 
Let $\Ga$ be the set of simple roots. For a fixed nondegenerate invariant
inner product $(.,.)$ on $\g,$ $\{ x_\be \, | \, \be \in \pm \De_+ \}$
will denote a set of root vectors $x_{\be} \in \g^{\be},$ 
such that 
\beq
(x_\be, x_{-\be}) = 1, \; \forall \be \in \De_+.
\label{norm}
\eeq

\bde{1.BDtriple}A Belavin--Drinfeld triple is a triple 
$(\Ga_1, \Ga_2, \tau)$ where $\Ga_1,$ $\Ga_2$ $\sub \Ga$ and
$\tau : \Ga_1 \ra \Ga_2$ is such that

1. $(\tau(\al), \tau(\be)) = (\al, \be)$ for all simple roots
$\al$ and $\be$ from the set $\Ga_1,$ 

2. For any $\ga \in \Ga_1$ there exists a positive integer $n$ for which
$\tau^n (\ga) \in \Ga_2 \backslash \Ga_1$ (nilpotency condition).
\ede
Given a Belavin--Drinfeld triple $(\Ga_1,$ $\Ga_2$ $\sub \Ga),$
one defines a partial ordering on $\De_+$ by $ \al < \be$ if
$\al \in \Ga_1,$ $\be \in \Ga_2,$ and $\be = \tau^n(\al)$
for an integer $n.$

\bth{1.BDclass}{\em{(}}Belavin--Drinfeld{\em{\/)}} Any factorizable 
Lie bialgebra structure on a complex reductive Lie algebra $\g$ has an
$r$-matrix
\beq
r = r_0 + \sum_{\al \in \De_+} x_{-\al} \o x_\al
        + \sum_{\al, \be \in \De_+, \\ \al < \be}
                               x_{-\al} \wedge x_\be
\label{BDr-matr}
\eeq
for an appropriate choice of a Cartan subalgebra $\h,$ 
an invariant bilinear form $(.,.)$ on $\g,$ root vectors 
$x_\be \in \g^\be$ normalized by \eqref{norm}, and a Belavin--Drinfeld 
triple $(\Ga_1, \Ga_2, \tau).$ 
The component $r_0 \in \h^{\o 2}$ of $r$ has to be such that
\beqa
&&(\tau \al \o 1) r_0 + (1 \o \al) r_0 =0, 
\; \forall \al \in \Ga_1,
\nn \\
&& r_0 + r_0^{21} = \Om_0.
\nn
\eeqa
Here $\Om_0$ is the component in $\h \o \h$ of the Casimir of $\g$
associated to the form $(.,.).$ 
\eth

For such an $r$-matrix the bilinear form $(.,.)_r$ on $\g$
induced by $r + r_{21}$ (see Sect.~2) is equal to the form $(.,.).$

The original result of Belavin and Drinfeld dealt with the case
of nonskewsymmetric $r$-matrices on complex simple Lie algebras.
It was later observed (see \cite{Ho}) that their proof
can be easily extended to the above generality.

It is also important to note that all Lie bialgebra structures
on a complex semisimple algebra $\g$ are quasitriangular
since $H^2( \g, \g)= 0.$ But this is not true for a reductive
algebra. For example when $\g$ is commutative any Lie algebra
structure on $\g^\st$ gives rise to Lie bialgebra structure on
$\g.$

Let us fix an $r$-matrix as in \thref{1.BDclass} associated to
a Belavin--Drinfeld triple
$(\Ga_1, \Ga_2, \tau).$ Denote with $\p_+$ the parabolic subalgebra of
$\g$ containing $\b_+$ and the root spaces $\g^{\al}$ for all negative
roots $\al$ that are linear combinations of simple roots from $\Ga_1.$
In a similar way define a parabolic subalgebra $\p_-$ containing
$\b_-$ and $\g^{\be}$ for all positive roots $\be$ that are
linear combinations of simple roots from $\Ga_2.$

With $\n_+,$ $\n_-$ we will denote the unipotent radicals
of $\p_+,$ $\p_-$ and with $\l_1,$ $\l_2$ their Levi factors. 

The Levi decompositions of $\p_\pm$ read:
\beq
\p_\pm = \l_{1,2} \lop \n_\pm. \label{1.lev1} 
\eeq
The derived subalgebras of $\l_1,$ $\l_2$ will be denoted by $\l'_1,$
$\l'_2$ and their centers by $\c_1,$ $\c_2,$ so
\beq
\l_i = \l'_i \op \c_i, \; i=1,2. \label{1.red} 
\eeq

The spans of first and second components of $r_0$ contain $\l'_1\cap \h$
and $\l'_2 \cap \h$ respectively. Denote by $\h_1$ and $\h_2$
their intersections with $\c_1$ and $\c_2.$ Finally define 
\beq
\h^\ort_i = \left[ \h_i \op (\l'_i \cap \h) \right]^\perp, \; i=1,2
\label{1.ortcompl}
\eeq
where the ortogonal complement is taken in $\h$ with respect to
the form $(.,.).$ One proves \cite{BD, ESh} that $\m_\pm$ are the
orthogonal complements to $\g_\pm$ which implies
\[
\h^\ort_i \sub \h_i, \; i=1, 2.
\]
Therefore $\g_\pm/ \m_\pm$ can be identified with 
$\l'_{1,2} \op (\h_+ / \h^\ort_+)$ and the Cayley
transform $\th$ can be considered as an orthogonal
Lie algebra
isomorphism between the last two algebras. 
One can choose complements $\a_i$ for $\h^\ort_i$ in $\h_i$ $(i=1,2),$
for which the map $\th$ can be lifted to a Lie algebra isomorphism
\beq
\th: \; \g'_1 \op \a_1 \ra \g'_2 \op \a_2,
\label{1.theta}
\eeq
preserving the form $(.,.)$ in such a way that $\th(\a_1) = \a_2.$

\bpr{1.BDalg} For a Belavin--Drinfeld triple $(\Ga_1, \Ga_2, \tau)$
the algebras $\g_\pm,$ $\m_\pm,$ are given by
\beqa
&&\m_\pm = \h^\ort_{1,2} \lop \n_\pm,
\label{1.m} \\
&&\g_\pm = (\l'_{1,2} \op \a_{1,2}) \lop \m_\pm, 
\label{1.g} 
\eeqa
and thus the quotients $\g_\pm / \m_\pm$ are 
\beq
\g_\pm/ \m_\pm \cong \l'_{1,2} \op \a_{1,2}.
\label{1.C} 
\eeq
The Cayley transform $\th: \l'_1 \op \a_1 \ra \l'_2 \op \a_2$ maps
$\l'_1$ to $\l'_2$ extending the map $\tau$ to root spaces 
and $\a_1$ to $\a_2$ preserving the form $(.,.)$ {\em{(}}where it is
completely determined by the term $r_0$ of the $r$-matrix $r${\em{)}}. 
\epr

The projections 
$\pr_\pm: \; \g_\pm \ra \l'_{1,2} \op \a_{1, 2} \cong \g_\pm / \m_\pm$ 
(recall eq. \eqref{1.m}) are just the projections on the first components
in \eqref{1.g}.

The decompositions \eqref{1.g} are similar to the Levi decompositions of
$\g_\pm.$ The Lie algebras $\m_\pm$ stay ``between'' the nilpotent
radicals $\n_\pm$ of $\g_\pm$ and their radicals $\c_{1,2} \lop \n_\pm.$

\bex{1.fullh} Consider the special case when the $r$-matrix $r$ is
such that 
\[
\g_+ \sup \h.
\]
This means that the Cartan subalgebra $\h$ sits entirely inside the
subalgebra of $\g$ spanned by the first components of $r$ and is
equivalent to the condition $\h^\ort_i = 0.$ Then $\g_\pm$ are full
parabolic subalgebras
\[
\g_\pm = \p_\pm
\]
and their ideals $\m_\pm$ are equal to the unipotent radicals
$\n_\pm$ of $\p_\pm$
\[
\m_\pm = \n_\pm.  
\]
The quotient algebras $\g_\pm /\m_\pm$ are isomorphic to the Levi factors
$\l_{1,2}$ of $\p_\pm:$
\[
\g_\pm /\m_\pm \cong \l_{1,2}.
\]
So $\a_i$ can be naturally taken $\c_i$ $(i=1, 2)$ and the Cayley
transform $\th : \: \l_1 \ra \l_2$ is an orthogonal Lie algebra 
isomorphism between the full Levi factors $\l_1$ and $\l_2.$ 
\eex

\bex{1.stand}The {\em{standard}} Lie bialgebra structure on
a complex semisimple Lie algebra $\g$ comes from the $r$-matrix
\[
r = \sum_{x_1} x_i \o x_i + \sum_{\al \in \De_+} x_{-\al} \o x_\al
\]
where $\{ x_i \}$ is an orthonormal basis of $\h$ with respect to the
Killing form on $\g.$ It is associated with the trivial Belavin--Drinfeld
triple for which $\Ga_1 = \Ga_2 = \emptyset.$ The algebras
$\g_\pm$ are equal to the Borel subalgebras $\b_\pm,$
and $\m_\pm$ are their unipotent radicals $\n_\pm.$
The quotients $\g_\pm/ \m_\pm$ are naturally isomorphic
to the Cartan subalgebra $\h$ $(\l'_i$ and $\h^\ort_i,$ $i=1,2$ are
trivial). Under this identification the map $\th: \h \ra \h$
is $\th = - \id_\h.$
\eex

\bex{CG} Let $\g = \s \l (n+1),$ $\h$ be its Cartan
subalgebra consisting of diagonal matrices, and 
$\{ \al_1, \dots, \al_n \}$ be the set of its simple roots.

The Cremmer--Gervais \cite{CG, Ho} Lie bialgebra structure
on it is associated to the Belavin--Drinfeld triple
$\Ga_1 = \{ \al_1, \dots, \al_{n-1} \},$
$\Ga_2 = \{ \al_2, \dots, \al_n \},$ $\tau(\al_j) = \al_{j+1},$
$j=1, \ldots, n-1.$ The algebras $\h^\ort_i$ are trivial and
\beqa
&&\l_1 = \{ \diag(A, - \tr A) | A \in \g \l (n) \} \cong \g \l(n),
\label{CGl1} \\
&&\l_2 = \{ \diag(- \tr A, A) | A \in \g \l (n) \} \cong \g \l(n)
\label{CGl2}
\eeqa
(in block notation). The $r_0$ component of the $r$-matrix is chosen 
in such a way that its Cayley transform $\th: \l_1 \ra \l_2$ is 
\[
\th( \diag(A, -\tr A)) = \diag(-\tr A, A).
\]
\eex
%%%%%%%%%%%%%%%%%%%%%%%%%%%%%%%%%%%%%%%%%%%%%%%%%%
\subsection{Related Poisson--Lie groups}
%%%%%%%%
For a complex reductive group $R$ we will denote its derived 
group by $R'$ and the identity component of its center by $Z(R)^\ci.$
Recall the standard fact
\beq
R = R' \: Z(R)^\ci
\label{reductive}
\eeq 
and that $R' \cap Z(R)^\ci$ is a finite abelian group.

In the rest of this paper $G$ will denote a connected  
complex reductive Lie group. We will fix a factorizable Lie bialgebra
structure on  $\g:= \Lie(G)$ associated to a Belavin--Drinfeld triple
$(\Ga_1, \Ga_2, \tau)$ and equip $G$ with the corresponding
Sklyanin Poisson--Lie structure (see \eqref{sklyanin}).
 
With $G_+$ and $G_-$ we will denote the connected
subgroups of $G$ with tangent Lie bialgebras $\g_+$ and $\g_-.$ 
They are Poisson--Lie subgroups of $G$ since $\g_\pm$ are
Lie sub-bialgebras of $\g.$ Let $M_\pm$ be the subgroups of $G_+$ with
$\Lie(M_\pm)= \m_\pm.$ 

Assuming the notation from Subsect.~3.1 for $\Lie(G)=\g,$
we will give explicit formulas for $M_\pm$ and $G_\pm,$ 
similar to \eqref{1.m} and \eqref{1.g}.
The maximal
torus of $G$ corresponding to the fixed Cartan subalgebra $\h$
of $\g$ will be denoted by $H.$ Let $B_\pm$ be the Borel subgroups
of $G$ relative to $H$ with $\Lie(B_\pm) = \b_\pm.$ Denote with $P_\pm$
the parabolic subgroups of $G$ with Lie algebras $\p_\pm.$ Let $L_{1,2}$
be their Levi factors and $N_\pm$ be their unipotent radicals. The
Levi decompositions of $P_\pm$ are
\[
P_\pm = L_{1,2} \ltimes N_\pm.
\]

The connected subgroups of $Z(L_i)^\ci$ with Lie algebras
$\a_i$ and $\h^\ort_i$ will be denoted by $A_i$ and $H^\ort_i.$    

The group analogs of formulas \eqref{1.m} and \eqref{1.g}
are
\beqa
&&M_\pm = H^\ort_{1,2} \ltimes N_\pm,
\label{Mpm} \\
&&G_\pm = (L'_{1,2} \: A_{1,2}) \: M_\pm.
\label{Gpm}
\eeqa
The intersections $L'_i \cap A_i \sub Z(L'_i)$ are finite groups, but they
are nontrivial in general. The groups $L'_{1,2} \: A_{1,2}$ normalize
$M_\pm.$ The intersections
\[
(L'_{1,2} \: A_{1,2}) \cap M_\pm =
(Z(L'_{1,2}) \: A_{1,2}) \cap H^\ort_{1,2}
\]
are discrete subgroups of $Z(L'_{1,2} A_{1,2}).$ Generally
they are also nontrivial.

The factor groups $G_\pm / M_\pm$ are
\beq
G_\pm / M_\pm \cong 
(L'_{1,2} \: A_{1,2}) / 
(Z(L'_{1,2})\, A_{1,2} \cap H^\ort_{1,2}) ).
\label{F12} 
\eeq
Unfortunately the Cayley transform 
$\th: \l'_1 \op \a_1 \ra \l'_2 \op \a_2$
cannot be lifted in general to an isomorphism between the two groups
in \eqref{F12}. One of the reasons for this is that 
it migth not map the kernel of the exponential map
for $L'_1 A_1$ to the one of the exponential
map for $L'_2 A_2.$
Nevertheless there exist two discrete subgroups
$\La_i$ of $Z(L'_i \: A_i)$ containing $(L'_i A_i) \cap H^\ort_i$ such
that the factor groups $(L'_i \: A_i)/ \La_i$ are isomorphic 
(as abstract groups
only). We will consider the minimal such groups $\La_i.$ Explicitly they
are given by 
\[
\La_i = 
\Big( \exp_i \circ \, \th^{-\e_i} \circ \exp_{3-i}^{-1}
  ( Z(L'_{3-i})A_{3-i} \cap H^\ort_{3-i}) \Big)
  (Z(L'_i)A_i \cap H^\ort_i)
\]
where $\e_i =(-1)^i,$ $i=1,2.$ The corresponding isomorphism 
will be denoted by $\Th:$
\beq
\Th: (L'_1 \: A_1)/ \La_1 \ra (L'_2 \: A_2)/ \La_2.
\label{Theta2}
\eeq
Clearly the groups $\La_i$ consist of semisimple elements
of $L'_i A_i.$ Thus the restrictions of the projections
$L'_i A_i \ra (L'_i A_i)/\La_i$ to any unipotent subgroup
of $L'_i A_i$ are one to one. 

\bre{LaLa}Note that the groups $\La_i$ might have in general some
accumulation points and then $(L'_i \: A_i)/ \La_i$ will not have
structure of differential manifolds. In Sect.~3.2.1 and Sect.~4 we will
assume that the $r$-matrix, we start with, is such that this does not
happen.   
In Sect.~5 the full generality will be assumed. It brings the
inconvenience that $G^r$ is not a closed subgroup of $G \times G,$
but this will not be essential for the proofs there.
\ere
Let
\[
\Pr_\pm: G_\pm \ra (L_{1,2} \: A_{1,2}) / \La_{1,2}
\]
denote the compositions of the projection maps
\[
G_\pm \ra G_\pm / M_\pm \cong 
(L'_{1,2} \: A_{1,2}) /
((L'_{1,2} \: A_{1,2}) \cap H^\ort_{1,2}) )
\ra (L'_{1,2} \: A_{1,2})/ \La_{1,2}.
\]
Their differentials at the origin are equal to the maps
$\pr_\pm: \g_\pm \ra \g_\pm / \m_\pm$ introduced in Sect.~2.

\bre{fullHgr} If the Poisson structure on $G$ comes from 
an $r$-matrix $r$ satisfying $\g_+ \sup \h$ (see \exref{1.fullh}),
then the
groups $G_\pm$ are the full parabolic subgroups $P_\pm$ of $G,$
$M_\pm$ are their unipotent radicals $N_\pm,$ and 
$L'_{1,2} A_{1,2}$ are their Levi factors $L_{1,2}.$
$(H^\ort_{1,2}$ are trivial.) The groups $\La_{1,2}$ are generally
nontrivial even if $P_\pm = B_\pm$ and $L_i = H.$  
\ere

\bex{2.stand}The standard Poisson structure on a semisimple
Lie group $G$ is the integration of the standard Lie bialgebra structure
on $\g = \Lie(G)$ reviewed in \exref{1.stand}.
In this case
$G_\pm = B_\pm,$ $M_\pm = N_\pm$ where $N_\pm$ are the unipotent
radicals of the Borel subgroups $B_\pm.$
The groups $L'_{1,2},$ $H_{1,2}^\ort$ are trivial and $A_{1,2}= H.$
There is no need to introduce additional groups $\La_i$ and
the Cayley transform $\Th: H \ra H$ acts as
$\Th(h) = h^{-1}.$
\eex
\bex{CGgroup} The groups $L_i$ for the Cremmer--Gervais
Poisson structure on $SL(n+1)$ (see \exref{CG}) are
\beqa
&& L_1 = \{ \diag(A, (\det A)^{-1} ) | A \in GL(n) \} \cong GL(n), 
\nn \\
&& L_2 = \{ \diag( (\det A)^{-1}, A ) | A \in GL(n) \} \cong GL(n).
\nn
\eeqa
$G_\pm$ are the parabolic subgroups of $SL(n+1)$ containing $B_\pm$
with Levi factors $L_{1, 2}$ respectively.
The map $\th: \l_1 \ra \l_2$ lifts to $\Th: L_1 \ra L_2:$
\[
\Th( \diag(A, (\det A)^{-1} ) ) = \diag( (\det A)^{-1}, A),
\]
so $\La_i$ are trivial.
\eex
%%%%%%%%%
\subsubsection{The double $\D(G_+)$}
%%%%%%%
Combining the results from Subsect.~2.2 and eq. \eqref{F12},
we see that we can choose as a double Poisson--Lie group
$\D(G_+)$ of $G_+$ the following group
\beq
\D(G_+) := G \times ((L'_1 \: A_1)/ \La_1).
\label{dobll}
\eeq
The groups $G_+$ and $G_-$ are embedded in it by
\beqa
&&i_+(x) = (x, \Pr_+ (x)), \; x \in G_+,
\label{G+embed} \\
&&i_-(y) = (y, \Th^{-1} \circ \Pr_- (x)), \; y \in G_-
\label{G-embed}
\eeqa
(cf. \prref{1.embed}).
The group \eqref{dobll} is also a double of $G_-,$ but the corresponding
Poisson bivector field is opposite to the one for $\D(G_+)$ 
since $\g_- \cong \g_+^{\st \opp}.$
%%%%%%%%
\subsubsection{The double $\D(G)$}
%%%%%%%%
The results from Subsect.~2.1 imply that the double
Poisson--Lie group $\D(G)$ of $G$ can be chosen as 
\[
\D(G) := G \times G.
\]
$G$ is embedded in it diagonally:
\[
d: \: G \ra G \times G, \quad d(g) = (g, g), \; \mbox{for} \; g \in G.
\]
Its image will be denoted by $d(G).$ The connected
subgroup of $G \times G$ with Lie algebra $\g^r$ will
be denoted by $G^r$ (see Propositions \ref{pRS} and \ref{p1.g^r}). 
It plays the role of opposite dual Poisson--Lie group
of $G$ (see \thref{1.STS}) and is explicitly given by
\beq
G^r:= \left\{ (g_+, g_-) \; | \;
              g_+ \in G_+, \; g_- \in G_-, \; \mbox{and}
              \; \Th \circ \Pr_1(g_+) = \Pr_-(g_-)
     \right\}.
\label{Gr}
\eeq
Because of \prref{1.g^r} we need only to prove that the group in
\eqref{Gr} is connected. Consider its identity
component $(G^r)^\ci$ and define 
$\check{\La}_i \sub L'_i A_i$ by
\[
\check{\La}_i = \{ x \in L'_i A_i | (x, e) \in 
(G^r)^\ci \}, \; i=1,2.
\]
The following map is an isomorphism from
$(L'_1 A_1)/\check{\La}_1$ to $(L'_2 A_2)/\check{\La}_2:$
\[
l_1 \check{\La}_1 \mapsto l_2 \check{\La}_2, \;
\mbox{if} \; l_i \in L'_i A_i, \; \mbox{and} \;
(l_1, l_2) \in (G^r)^\ci.
\]
It is correctly defined since $L'_1 A_1$ is connected.
Besides this $\check{\La}_i \sup (L'_i A_i) \cap H^\ort_i$ because
$H^\ort_1 \times e$ and $e \times H^\ort_2$ are subgroups of
$(G^r)^\ci.$ From the minimality property of $\La_i$ we get
$\La_i \equiv \check{\La}_i$ which implies
$G^r = (G^r)^\ci.$
%%%%%%%%%%%%%%%%%%%%%%%%%%%%%%%%%%%%%%%%%%%%
\subsection{Induction on factorizable complex
reductive Lie bialgebras and \\ Poisson--Lie groups}
As in Subsect.~3.1 we will assume that $\g$ is a complex reductive
Lie bialgebra with an $r$ matrix as in \eqref{BDr-matr} associated
to the Belavin--Drinfeld triple $(\Ga_1, \Ga_2, \tau).$
Recall that the nilpotency condition for the map $\tau$ states that
for each $\ga \in\Ga_1$ there exists an integer $n(\ga),$
such that $\tau^{n(\ga)} \in \Ga_2 \backslash \Ga_1.$ Denote
\[
\ord(\tau)= \max_{\ga \in \Ga_1} n(\ga).
\]
We will canonically construct a chain of complex reductive
factorizable Lie
bialgebras $\g^{(k)}$ $(k = 0, \ldots, \ord(\tau))$ associated to   
Belavin--Drinfeld triples $(\Ga^{(k)}_1, \Ga^{(k)}_2, \tau^{(k)})$  
such that $\g^{(0)} \equiv \g$ as Lie bialgebras and
\[
\ord(\tau^{(k)})= \ord(\tau) -k.
\]
Recall that $\D(\g_+) \cong \g \op (\l'_1 \op \a_1)$
and $\g_+$ and $\g_-$ are embedded in it via the maps $i_\pm$
\eqref{1.incl+}--\eqref{1.incl-}. Let $\dmf$ be the subalgebra
of $\g \op (\l'_1 \op \a_1)$ defined by
\beq
\dmf = 
\l_2 \op ((\l'_1 \op  \a_1) \cap \l_2) =
\l_2 \op ((\l'_1 \cap \l_2) \op \a_1).
\label{falg}
\eeq
It is clear that the restriction of the standard bilinear form of
$\D(\g_+)$ (see the proof of \prref{1.embed}) to $\dmf$ is nondegenerate.
Define
\beqa
&& \g^{(1)}_+= ((\l'_1 \cap \l_2)  \op \a_1)
               \lop (\h_1^\ort \lop \n^{(1)}_+) \sub \g_+,
\label{g1+} \\
&& \g^{(1)}_-= (\th(\l'_1 \cap \l_2) \op \a_2)
               \lop (\h^\ort_2 \lop \n^{(1)}_-) \sub \g_-
\label{g1-}
\eeqa
where $\n^{(1)}_+,$ $\n^{(1)}_-$ are the nilpotent subalgebras of
$\l'_2$ spanned by the root spaces of positive / negative roots of
$\l'_2$ that are not roots of $\l'_1 \cap \l'_2,$
$\th(\l'_1 \cap \l'_2).$

The images $i_\pm(\g^{(1)}_\pm)$ of $\g^{(1)}_\pm$ under the
inclusions $\i_\pm$ (see \eqref{1.incl+}--\eqref{1.incl-}) are isotropic
subalgebras of $\dmf$ with trivial intersection since the same is
true for the full images $\i_\pm(\g_\pm)$ as subalgebras of $\D(\g_+).$
Taking into account $\dim \dmf = 2 \dim \g^{(1)}_\pm$ we get that
$(\dmf, i_+(\g^{(1)}_+), i_-(\g^{(1)}_-))$ is a Manin triple.
This induces Lie bialgebra structures on $\g^{(1)}_\pm$ and $\dmf$
such that $\dmf \cong \D(\g^{(1)}_\pm).$ One easily checks that the
second term of $\dmf$ in \eqref{falg} is a Lie bialgebra ideal of 
$\dmf.$
The factor of $\dmf$ by it is isomorphic to $\l_2$ as a Lie algebra
which induces a Lie bialgebra structure on $\l_2.$ Set
\[
\g^{(1)} := \l_2.
\]
Define the subsets $\Ga^{(1)}_1$ and $\Ga^{(1)}_2$ of $\Ga_2$ by
\beqa
&&\Ga^{(1)}_1 = \Ga_1 \cap \Ga_2,
\nn \\
&&\Ga^{(1)}_2 = \tau(\Ga_1 \cap \Ga_2),
\nn
\eeqa
and let
\beq
\tau^{(1)} := \tau|_{\Ga^{(1)}_1}:
\: \Ga^{(1)}_1 \ra \Ga^{(1)}_2.
\label{1tau}
\eeq
The definition of the Lie bialgebra structure on $\g^{(1)}$ as
a quotient of the double $\dmf$ \eqref{falg} implies that
it comes from the $r$-matrix
\beq
r^{(1)}= r_0 + \sum_{\al \in \De_{2+}} x_{-\al} \o x_\al
+ \sum_{\al, \be \in \De_{2+}, \al < \be}
x_{-\al} \wedge x_{\be}
\label{r1matr}
\eeq
where $\De_{2+}$ denotes the set of positive roots of $\l_2$ and
$<$ is the restriction of the initial partial ordering of $\De_+$
to $\De_{2+}$ (see \thref{1.BDclass}).
The Belavin--Drinfeld triple associated to the
$r$-matrix $r^{(1)}$ is $(\Ga^{(1)}_1, \Ga^{(1)}_1, \tau^{(1)}).$
The images of the corresponding homomorphisms
$r^{(1)}_\pm: \: (\g^{(1)})^\st \ra \g^{(1)}$
are exactly the Lie bialgebras $\g^{(1)}_+$ and $\g^{(1)}_-.$ From
\eqref{1tau} it is clear that
\[
\ord(\tau^{(1)}) = \ord(\tau) -1.
\]

Repeatedly applying this construction, one builds a sequence 
$\g^{(0)} = \g,$ $\g^{(1)},$ $...,$ $\g^{(ord \tau)}$
of reductive factorizable Lie bialgebras with the stated property.

The Lie bialgebra $\g^{(ord \tau)}$ obtained at the end
is the Cartan subalgebra $\h$ of $\g$ with trivial Lie cobracket
(since it is an abelian coboundary bialgebra). Besides this
$\g^{(\ord \tau)}_\pm$ and $\D(\g^{(\ord \tau)})$ are
given by
\beqa
&&\g^{(\ord \tau)}_\pm = \g_\pm \cap \h,
\nn \\
&&\D(\g^{(\ord \tau)}) =
\h \op ((\l'_1 \op \a_1)\cap \h) =
\h \op \left( (\l'_1 \cap \h) \op \a_1 \right).
\nn
\eeqa

On the group level, one constructs a chain of connected
complex reductive Poisson--Lie groups $G^{(0)}= G,$
$G^{(1)},$ $\ldots,$ $G^{(\ord \tau)}=H.$ As a Lie group
$G^{(1)}$ is just $L_2$ and $G^{(1)}_\pm$ are given by
\beq
G^{(1)}_\pm= \left( L^{(1)}_{1,2} \: A^{(1)}_{1,2} \right) \,
\left(H_{1,2}^\ort \ltimes N^{(1)}_\pm \right),
\label{G+-1}
\eeq
where $L^{(1)}_1,$ $L^{(1)}_2,$ $A^{(1)}_1,$ and $A^{(1)}_2$
denote the connected subgroups of $G$ with Lie algebras
$\l'_1 \cap \l'_2,$ $\th(\l'_1 \cap \l'_2),$
$(\l'_1 \cap \c_2) \op \a_1,$ and $\th(\l'_1 \cap \c_2) \op \a_2.$
$N^{(1)}_\pm$ are the unipotent subgroups of $G$ with
Lie algebras $\n^{(1)}_\pm.$

The corresponding groups $M^{(1)}_\pm \sub G^{(1)}_\pm$
(see Subsect.~3.2) are the second factors in these formulas.

On Lie algebra level the Cayley transform
$\th^{(1)}: \Lie(L^{(1)}_1 \, A^{(1)}_1) \ra \Lie(L^{(1)}_2 \, A^{(1)}_2)$
is the restriction of $\th$ to the subspace
$\Lie(L^{(1)}_1 \, A^{(1)}_1)$ $=$ $(\l'_1 \op \a_1) \cap \l_2$
of $\l'_1 \op \a_1.$ Since $\La_i$ are subgroups of
$Z(L'_i \, A_i) \sub H$ they sit inside 
$L^{(1)}_i \, A^{(1)}_i$ and $\th^{(1)}$ can be lifted
to an isomorphism
\[
\Th^{(1)}: \:
(L^{(1)}_1 \: A^{(1)}_1)/ \La_1 \ra
(L^{(1)}_2 \: A^{(1)}_2)/ \La_2.
\]
The double $\D(G^{1}_\pm)$ can be taken as 
\beq
\D(G^{(1)}_\pm): =G^{(1)} \times
(L^{(1)}_1 \: A^{(1)}_1 / \La_1).
\eeq
It is naturally a subgroup of $G \times (L_1 \, A_1 /\La_1)$
and the embeddings
\[
i^{(1)}_\pm: \: G^{(1)}_\pm \ra \D(G^{(1)})
\]
are simply the restrictions of the embeddings $i_\pm$
(see \eqref{G+embed}--\eqref{G-embed}) to $G^{(1)}_\pm.$
 
The last object $G^{( \ord \tau)}$ is
\[
G^{(\ord \tau)} = H
\]
and
\[
G^{(\ord \tau)}_\pm = (L'_{1,2}\cap H) \: A_{1,2} \: H^\ort_{1,2}.
\]

At each step of the described procedure
there is a second canonical way of constructing a Lie bialgebra
$\g^{k+1}$ from $\g^k,$ having the stated property.
Restricting to the first step $(k=0),$ define $\g^{(1)} = \l_1,$
\beqa
&& \g^{(1)}_+ = (\th^{-1}(\l'_2 \cap \l_1) \op \a_1) \lop 
                (\h_1^\ort \lop \th^{-1}(\n^{(1)}_1)),
\nn \\
&& \g^{(1)}_- = (\l'_2 \cap \l_1 \op \a_2) \lop 
                (\h_2^\ort \lop \th^{-1}(\n^{(1)}_1)),
\nn
\eeqa
and induce Lie bialgebra structure on $\g^{(1)}$ using the subalgebra
$\l_1 \op ((\l'_2 \op  \a_2) \cap \l_1)$ of $\D(\g_+).$
This structure comes from an $r$-matrix as in \eqref{r1matr}
with sums over the roots of $\l'_1$ insted of $\l'_2.$
The two choices for Lie bialgebras $\g^{(1)}$ are isomorphic
if $\h_\pm^\ort$ are trivial, but seem to be slightly different in general. 
%%%%%%%%%%%%%%%%%%%%%%%%%%%%%%%%%%%%%%%%%%%%%%%%%%%%%%%%%%%%%%%%%%%%%%%%%%%%
\sectionnew{Symplectic leaves in $G_-$}
%%%%%%%%%%%%%%%%%%%%%%%%%%%%%%%%%%%%%%%%%%%%%%%%%%%%%%%%%%%%%%%%%%%%%%%%%%%%
Subsect.~4.1 is devoted to the description of the set of
$i_+(G_+)$ double cosets of $\D(G_+),$ while Subsect.~4.2
studies the discrete subgroup 
$i_-(G_-) \cap i_+(G_+)$ of $\D(G_+).$ In the last Subsection these 
results are combined with \thref{int} to describe the set of 
symplectic leaves of $G_-.$
%%%%%%%%%
\subsection{Double cosets of $i_+(G_+)$ in $\D(G_+)$}
%%%%%%
Denote the double coset of the element
\[
(g, \til{l}) \in G \times \left((L'_1 A_1) / \La_1\right) \cong
\D(G_+)
\]
by $[(g, \til{l})]:$ 
\[
[(g,\til{l})] = i_+(G_+) (g,\til{l}) i_+(G_+).
\]
Since $L'_1 \, A_1 \sub G_+,$ for all $g \in G,$ 
$\til{l} \in L'_1 \, A_1$
\beq  
[(g, \pr_+(l) \til{l})]= [(l^{-1} g, \til{l})] \; \mbox{and} \;
[(g, \til{l} \pr_+(l))]= [(g l^{-1}, \til{l})], 
\; \forall l \in L'_1 A_1.
\label{compute}
\eeq

The following Proposition illustrates the problem of classifying
double cosets.

\bpr{G+1}The $i_+(G_+)$ double cosets of $\D(G_+)$ are in one to one
correspondence with the conjugation orbits of $G_+$ in 
$G / (\La_1 \, M_+).$
\epr
Recall that $M_+$ is a normal subgroup of $G_+$ and $\La_1 \in Z(L_1),$ so
the conjugation action of $G_+$ on $G$ indeed can be pushed down to the
the homogeneous space $G /\La_1 M_+.$ Note also that this action can
equivalently be described as the action of $(L'_1 \, A_1) \, M_+$ on
$G/(\La_1 \, M_+)$ where the first factor acts by conjugation and the
second one by left multiplication. \\
{\em{Proof of \prref{G+1}.}} By a straightforward computation for any 
$g \in G$  
\beqa
[(g, e)]
&=& \left\{ (m' l' g l'' m'',
             \pr_+(l') \pr_+(l'')) |  \;
             l', l'' \in L'_1 A_1, \;
             m', m'' \in M_+
    \right\}
\nn \\
&=& \left\{ (m' l' g (l')^{-1} m'' k, \pr_+(k)) |  \;
             l', k   \in L'_1 A_1, \;
             m', m'' \in M_+
    \right\}.
\label{comp1}
\eeqa
The second expression is obtained from the first by letting
$l' l'' = k \in L'_1 \, A_1$ and expressing $l''$ in terms of $l'$ and
$k.$ We also used the fact that $M_+$ is a normal subgroup of $G_+.$

To the  $G_+$ conjugation orbit of $g(\La_1 \, M_+)$ we associate the
double coset $[(g, e)].$ Eq. \eqref{comp1} shows that 
this map is defined correctly. It is onto since each double coset of
$\D(G_+)$ is equal to a coset $[(g, e)]$ for some $g \in G,$ because of
\eqref{compute}. To show that it is injective assume that $[(g_1, e)] =
[(g_2, e)].$
If 
\[
(m'_1 l'_1 g_1 (l'_1)^{-1} m''_1 k_1, \pr_+(k_1))=
(m'_2 l'_2 g_2 (l'_2)^{-1} m''_2 k_2, \pr_+(k_2))
\]
is an elment of this double coset (the notation 
is as in \eqref{comp1}) then $k_1 \, (k_2)^{-1} \in \La_1.$
The first components give immediately that
the $G_+$ orbits of $g_1 \, \La_1 M_+$ and $g_2 \, \La_1 M_+$
coincide.
\hfill \qed

In the special case of an $r$-matrix $r$ for which
$\g_+ \sup \h$ (see \reref{fullHgr}), \prref{G+1} is about the conjugation
action of the parabolic subgroup $P_+$ of $G$ on $G/(\La_1 N_+),$
where $N_+$ is its unipotent radical. 

\subsubsection{Relation to double cosets in $G$ and the Bruhat Lemma}
Let
\beq
\Pi: \; \D(G_+) \cong G \times ((L'_1 \, A_1) / \La_1) \ra G.
\label{defPi}
\eeq
be the projection map on the first component of $\D(G_+).$
It is clear that
\beq
\Pi([(g, \til{l})]) = G_+ g G_+, \; 
\forall g \in G, \; \til{l} \in L'_1 \, A_1
\label{projPi}
\eeq
which connects $i_+(G_+)$ double cosets in $\D(G_+)$ with 
$G_+$ double cosets in $G.$

The group $G_+$ differs from the parabolic subgroup $P_+$ of $G$ just on
a part of the maximal torus $H$ (see \eqref{Gpm}). This makes the description
of the $G_+$ double cosets of $G$ an easy consequence from the one
for $P_+$ double cosets. The latter is given by the
Bruhat Lemma.

For any element $w \in W$ we will denote with $\dot{w}$ a representative
of it in the normalizer $N(H)$ of $H$ in $G.$ The conjugation and
the adjoint actions of $\dot{w}$ on $G$ and $\g$ will be denoted
by $\Ad_{\dot{w}}.$

\ble{Bruhat} {\em{(}}Bruhat{\em{\/)}} 
Let $P_1$ and $P_2$ be two parabolic subgroups
of $G,$ each of  which contains one of the Borel subgroups 
$B_+$ or $B_-.$ If $W_{1,2}$ denote the Weyl groups of their Levi
factors $L_{1,2}$ {\em{(}}naturally embedded in $W${\em{\/)}} then
\beq
G = \sqcup P_1 \dot{w} P_2,
\label{bruh}
\eeq
with $w$ running over a set of representatives in $W$ of all
cosets from $W_1 \backslash W / W_2.$
\ele
After commuting $\dot{w}$ with the part $Z(L_2)^\ci$ of the
torus $H \sub P_2,$ one gets
\[
G = \sqcup P_1 \dot{w} (L'_2 \ltimes N_2)
\]
where $N_2$ is the unipotent radical of $P_2$ and the disjoint union
is over the same subset of $W$ as in \eqref{bruh}.

The standard length function on $W$ will be denoted by $l(.).$
In the following Lemma we collect some
facts for minimal length representatives of
double cosets of the Weyl group $W$ which 
will be needed later. 

\ble{MinLe}In the notation of \leref{Bruhat}:

1. Each double coset from $W_1 \backslash W / W_2$ has a unique
representative $w \in W$ of minimal length. It can be equivalently
described as the representative $w \in W$ of the coset for which:

$(\st)$ For all positive roots $\al$ of $\Lie(L_1)$ and $\be$ of
$\Lie(L_2),$
$w^{-1}(\al)$ and $w(\be)$ are positive roots of $\g.$ \\
For any such $w \in W,$ $\l^w = \Lie(L_1) \cap \dot{w}^{-1}(\Lie(L_2))$
is a reductive subalgebra of $\g$ generated by $\h$ and root spaces
of simple roots of $\g.$ Its Weyl group will be denoted by $W^w.$

2. Any element $u \in W$ can be uniquely represented as
\beq
u = w_1 w w_2
\label{prodd}
\eeq
with $w,$ $w_1$ being minimal length representatives in $W,$ $W_1$ 
of cosets from $W_1 \backslash W / W_2,$ $W_1 / W^w$ and $w_2 \in W_2.$ 
Its length is given by
\beq
l(u) = l(w_1) + l(w) + l(w_2).
\label{sumlength}
\eeq

3. Let $P_3$ be a parabolic subgroup of $G$ such that 
$B_\pm \sub P_3 \sub P_1$ and $W_3$ be the Weyl group of its Levi 
factor. Then the set of minimal length representatives in $W$ for 
the cosets from
$W_3 \backslash W / W_2$ consists of the elements $u \in W$
whose representation from part 2 has  
$w_2 = \id$ and $w_1$ -- minimal length representative 
in $W$ of a coset from $W_3 \backslash W_1 / W^w.$
\ele
\proof In the case of $P_1=B_\pm$ (i.e. $W_1$ trivial) the Lemma is
proved in \cite{H}, Proposition~1.10. We will use this fact and will
refer to it as the single coset version of the Lemma.

Parts 1 and 2: Fix a coset from $W_1 \backslash W / W_2.$ Chose a
representative $w$ of it in $W$ having the property $(\st).$ Any minimal
length representative $w_0$ of the double coset has this property since
otherwise $l(s_1 w_0) < l(w_0)$ or $l(w_0 s_2) < l(w_0)$ for some simple
reflections $s_i \in W_i,$ $i=1,2.$

It is clear that any element $u \in W_1 w W_2$ can be represented in the
form \eqref{prodd} with $w_i \in W_i$ as in part 2 of the
Lemma. We will show the uniqueness of this representation 
and prove \eqref{sumlength}. The latter implies the 
statement of part 1.

Since $w_1 \in W_1$ is the minimal length representative of a coset from  
$W_1 / W^w,$ it has the property that $w_1(\al)$ is a positive root
of $\Lie(L_1)$ for any positive root $\al$ of $\l^w.$ This combined with
the definition of $w$ implies that $w_1 w$ is a minimal length
representative in $W$ of the coset $w_1 w W_2.$
(Here we also use the fact that $w_1$ maps positive roots of $\g$ that are
not roots of $\Lie(L_1)$ into positive roots of $\g.)$ The single coset
version of the Lemma applied to the coset $w_1 w W_2$ implies that
$w_2$ is determined uniquely from $u$ and therefore so is $w_1.$

The single coset version of \eqref{sumlength} for the cosets
$w_1 w W_2 \in W / W_2$ and $w_1 W^w \in W_1 /W^w$ with minimal
representatives $w_1 w$ and $w_1$ gives
\beqa
&& l(w_1 w w_2) = l(w_1 w) + l(w_2)
\nn \\
&& l(w_1 w) = l(w_1) +l(w)
\nn
\eeqa
which together imply \eqref{sumlength}.

Part 3: Part 1 implies that the elements of $W$ listed in part 3 are
minimal length representatives for cosets from $W_1 \backslash W / W_2.$
In the opposite direction, any minimal representative 
$u = w_1 w w_2 \in W$ written as in part 2 should have $w_2 = \id$ and
$w_1(\al)$ should be positive root of $\Lie(L_1)$ for any positive root
$\al$ of $\l^w.$ This implies the stated property of $w_2.$
\hfill
\qed

For the given $r$-matrix denote with $W_{1,2}$ the Weyl groups of the
Levi factors $L_{1,2}$ of $P_\pm.$ The Bruhat Lemma and \eqref{projPi} 
imply that any $i_+(G_+)$ double coset in $\D(G_+)$ is
equal to a coset $[(h \dot{w}, \til{l})]$ for some $h \in H,$ 
$\til{l} \in (L'_1 \, A_1) / \La_1$ and $w$ -- minimal length
representative in $W$ of a double coset from $W_1 \backslash W /W_1.$
Taking into account \eqref{compute} we obtain:

\ble{G_+interm}Each $G_+$ double coset of $\D(G)$ coincides with
one of the cosets $[(l \dot{w}, e)]$ for $w$ -- minimal representative
in $W$ of a coset from $W_1 \backslash W /W_1$ and $l \in L_1.$
\ele   
  
\subsubsection{Main result}
Let $V \sub W$ denotes the set of minimal length representatives
of the cosets from $W / W_1.$ According to \leref{MinLe}, each such coset 
has a unique representative in $V.$ 

For an element $v \in V$ denote by $\g^v$ the largest Lie subalgebra of
$\l_1$ that is stable under $\Ad_{\dot{v}}.$ It does not depend
on the choice of representative of $v$ in $N(H)$ and can also be
defined by:
\beq \g^v = \span \{ h, \; \g^\al \;
| \; \mbox{for the roots} \; \al \; \mbox{of} \; \l_1 \; \mbox{such that}
\; v^n (\al) \; \mbox{is a root of} \; \l_1 \; \forall \; \mbox{integer}
\; n \}.  
\label{lw} 
\eeq 
The minimality condition on $v$ guarantees that $\g^v$ is a reductive Lie
subalgebra of $\g,$ generated by the Cartan subalgebra $\h$ and root spaces
of simple roots of $\g.$ This follows from the fact that for any positive
root $\al$ of $\l_1,$ $v(\al)$ is a positive root of $\g$ (cf.
\leref{MinLe}). Denote with $G^v$ the reductive subgroup of $G$
with Lie algebra $\g^v$ and with $W^v$ its Weyl group. Clearly
$G^v$ is stable 
under $\Ad_{\dot{v}}.$

The definition of $\g^v$ combined with the minimality property of $v$
implies that for any positive root $\al$ of $\g^v,$ $v(\al)$ is again
a positive root of $\g^v,$ and thus so is $v^{-1}(\al).$ From part 1 of 
\leref{MinLe} it follows that any $v \in V$ is the minimal 
length representative of the double coset $W^v v W_1.$

The fact that $\Ad_{\dot{v}}$ preserves $G^v$ implies that it also
preserves its derived group $(G^v)'$ and the identity component of
its center $Z(G^v)^\ci.$ Thus we can define an action of $(G^v)'$
on itself by
\beq
TC^{\dot{v}}_{g}(f) = g^{-1} f \Ad_{\dot{v}}(g), \; g, \: f \in (G^v)'
\label{tc}
\eeq
which will be called {\em{twisted conjugation action.}}

The fact that $\La_1$ is a subgroup of $Z(L'_1 \, A_1) \sub H$
implies $\La_1 \sub G^v.$ Let $\La^v_1$ be the subgroup of
$(G^v)' \times Z(G^v)^\ci$ which is the inverse
image of $\La_1$ under the multiplication map
\[
(G^v)' \times Z(G^v)^\ci \ra G^v = (G^v)' Z(G^v)^\ci.
\]
Note that $\La^v_1 \sub Z((G^v)' \times Z(G^v)^\ci)$  and
that each element of $\La_1$ has finitely many
preimages in $\La^v_1$ (more precisely
$\# \, (G^v)' \cap Z(G^v)^\ci).$

Denote by $L^v$ the connected subgroup of $G^v$
with Lie algebra $\g^v \cap (\l'_1 \op \a_1).$
Let $\Cong^v(L^v)$ be the subgroup of the abelian
group $Z(G^v)^\ci$ defined by
\[
\Cong^v(L^v)= \{ h^{-1} \Ad_{\dot{v}}(h)
| \, h \in Z(L^v)^\ci \}.
\]
$H^\ort_1$ and $\Ad_{\dot{v}}(H^\ort_1)$ are also subgroups
of $Z(G^v)^\ci$ since by the definition of $H^\ort_1,$
\\ $H^\ort_1 \sub Z(L_1)^\ci.$ 

\bth{G+coset}The $i_+(G_+)$ double cosets of $\D(G_+)$ are classified
by a minimal length representative $v$ in $W$ of a coset from $W/ W_1,$
an orbit of the twisted conjugation action \eqref{tc} of $(G^v)'$ on
itself, and an element of the abelian group
\beq
Z(G^v)^\ci / ( \Cong^v(L^v) \, H^\ort_1 \, \Ad_{\dot{v}}(H^\ort_1) ).
\label{abelf}
\eeq
For a fixed  $v$ this data has to be taken modulo the multiplication
action of \\
$\La^v_1 \sub Z((G^v)' \times Z(G^v)^\ci).$

The dimension of the coset corresponding to the orbit through
$f \in (G^v)'$ is:
\beq
\dim \g_+ - \dim \h^\ort_1 +
\dim (\l'_1 \op \a_1) - \dim \l^v
+ l(v) + TC^{\dot{v}}_{(G^v)'}(f) +
\dim \Cong^v(L^v) H^\ort_1 \Ad_{\dot{v}}(H^\ort_1). 
\label{dim}
\eeq
{\em{(}}The described data corresponding to different choices of a
representative in $N(H)$ of a fixed element $v \in W$ can be
canonically identified{\em{.\/)}}
\eth
\proof In the previous Subsection we showed that any $i_+(G_+)$
double coset of $\D(G_+)$ is equal to a coset $[(l \dot{w}, e))]$
for some $l \in L_1$ and $w$ minimal length representative in $W$
of a coset from $W_1 \backslash W / W_1.$

Based on the elements $l$ and $\dot{w}$ we will
define a sequence of reductive Lie subalgebras $\g^k$ of $\g,$
elements $g^k \in \g^k,$ and $w^k \in W$ $(k \geq 1).$
It starts with $\g^1 = \l_1$ and $g^1=l,$ $w^1=w.$ 
The sequence will have the property  
that $\g^k$ and
\beq
\ol{\g}^k = \Ad_{(\dot{w}^{k-1} \ldots \dot{w}^1)} (\g^k)
\label{bargk}
\eeq
are reductive subalgebras
of $\g^{k-1} \sub \g$ containing the Cartan subalgebra $\h$ and
generated by root spaces corresponding to simple roots of $\g.$
If $W^k$ and $\ol{W}^k$ are their Weyl groups considered as
subgroups of $W^{k-1} \sub W,$ $w^k$ will have the property that it is
representative of minimal length of a coset from
$W^k \backslash W^{k-1} / \ol{W}^k.$

If $\Ad_{(\dot{w}^k \ldots \dot{w}^1)}$ does not preserve $\g^k$
define $\g^{k+1}$ by
\beq
\g^{k+1} = \g^k \cap \Ad_{(\dot{w}^k \ldots \dot{w}^1)}^{-1}(\g^k).
\label{gk+1}
\eeq
The minimality property of $w^j$ implies
that $w^k \ldots w^1$ is the minimal length representative 
of a coset from $W^k \backslash W / W^1$ (see part 3 of \leref{MinLe})
and therefore for any positive root
$\al$ of $\g^k,$ $(w^k \ldots w^1) (\al)$ and
$(w^k \ldots w^1)^{-1} (\al)$ are positive roots of $\g.$ From this we get
that both $\g^{k+1}$ and
$\ol{\g}^{k+1}=
\Ad_{(\dot{w}^k \ldots \dot{w}^1)}(\g^{k+1})$ are generated by
$\h$ and root spaces of simple roots of $\g.$

Let $\n_+^{k+1}$ and $\ol{\n}_+^{k+1}$ be the
two nilpotent subalgebras of $\g^k$ spanned by the root spaces
of positive roots of $\g^k$ that are not roots of
$\g^{k+1}$ and $\ol{\g}^{k+1}$ respectively.

The connected subgroups of $G$ with Lie algebras $\g^{k+1},$
$\ol{\g}^{k+1},$ $\n_+^{k+1},$ and $\ol{\n}_+^{k+1}$
will be denoted by $G^{k+1},$ $\ol{G}^{k+1},$
$N_+^{k+1},$ and $\ol{N}_+^{k+1}.$
They give rise to two parabolic subgroups
$G^{k+1} \ltimes N_+^{k+1}$ and 
$\ol{G}^{k+1} \ltimes \ol{N}^{k+1}_+$ of $G^k.$
To define recursively $g^{k+1}$ and $w^{k+1}$ from
$g^k$ and $w^k$ we apply the Bruhat Lemma for 
the above parabolic subgroups of $G^k.$
We get that $g^k \in G^k$ can be represented as
\beq
g^k = n' g' \dot{w}^{k+1} g'' n''
\label{wk+1}
\eeq
for some $w^{k+1} \in W^{k+1}$ that is the representative
of minimal length of a coset from \\ 
$W^{k+1} \backslash W^k / \ol{W}^{k+1}$
and some $g' \in G^{k+1},$ $g'' \in \ol{G}^{k+1},$
$n' \in N^{k+1}_+,$ $n'' \in \ol{N}^{k+1}_+.$
Eq. \eqref{wk+1} is the defining relation for $w_{k+1}$ 
and $g^{k+1}$ is defined by
\beq
g^{k+1} = \Big( \Ad_{(\dot{w}^k \ldots \dot{w}^1)}^{-1}(g'')^{-1} 
          \Big) g'.
\label{ggk+1}
\eeq
Clearly the described process continues for finitely many steps
$k=1, \ldots, K$ and at the end 
\beq
\Ad_{(\dot{w}^K \ldots \dot{w}^1)}(\g^K) = \ol{\g}^K = \g^K.
\label{olgg}
\eeq

Again due to part 3 of \leref{MinLe} the minimality properties
of $w^1,$ $\ldots,$ $w^K$ add up to
\beq
v:=(w^K \ldots w^1) \in V.
\label{v}
\eeq
Further we will use the representative
$\dot{v} = \dot{w}^K \ldots \dot{w}^1$ of $v$ in $N(H).$

Observe that
\[
\g^K = \g^v.
\]
Indeed, in view of \eqref{olgg} it is sufficient to prove
that $\g^v \sub \g^K.$ One proves by induction that $\g^v \sub \g^k$
for $k = 1, \ldots, K.$ The step from $k$ to $k+1$
follows from the definition \eqref{gk+1} of $\g^{k+1}$ and the fact
that $(w^K \ldots w^{k+2}) \in W^{k+1}.$ 

Thus we have $G^K = G^v = (G^v)' \, Z(G^v)^\ci$ and let
\beq
g^K = f^K c^K
\label{fc}
\eeq
for some $f^K \in (G^v)'$ and $c^K \in Z(G^v)^\ci$ defined
uniquely modulo multiplication of $f^K$ and division of
$c^K$ by an element of the finite group $(G^v)' \cap Z(G^v)^\ci.$
The statement from \thref{G+coset} about
the set of $i_+(G_+)$ double cosets follows from the following two
Lemmas.

\ble{lemma1}In the above notation
$[(l \dot{w}, e)] = [(g^K \dot{v}, e)]$ and two double cosets as the
one in the right hand side coincide if and only if they are related
to one and the same $v \in V,$ the twisted conjugation orbits
$TC^{\dot{v}}_{(G^v)'} (f^K)$ that correspond to them and the images
of the their elements $c^K$ in the factor group \eqref{abelf} are equal
modulo the multiplication action of $\La^v_1$ on
$(G^v)' \times Z(G^v)^\ci.$
\ele

For \thref{G+coset} we associate to the double coset
$[(l \dot{w}, e)] = [(g^K \dot{v}, e)]$ the $TC^{\dot{v}}$ orbit of
$f^K$ and the image of $c^K$ in the abelian factor group
\eqref{abelf} (recall eq. \eqref{fc}).
\ble{lemma2}Any element of $V$ can be obtained in a unique way as a
product $(w^K \ldots w^1)$ through the described iterative procedure. 
\ele

We will need some more notation: let $\n^k_-$
denote the nilpotent subalgebras of $\g^{k-1}$ spanned by
root spaces of negative roots of $\g^{k-1}$ that are not roots
of $\g^k.$ For convenience we set $\g^0=\g$ and define
$\n^k_\pm$ and $\ol{\n}^k_+$ for $k=1$ also.

Let also $\l^k= \g^k \cap (\l_1 \op \a_1)$ and denote by $L^k$ the
corresponding connected subgroup of $G.$ Note that
$L^K = L^v.$ \\ \hfill \\
{\em{Proof of \leref{lemma1}.}} Our strategy in dealing with the
double coset $[(l \dot{w}, e)]=[(g^1 \dot{w}^1, e)]$ is to commute
the right factor $i(G_+)$ with $\dot{w}=\dot{w}^1,$
splitting it according to how the result is related to
the left factor $i_+(G_+).$ Then we continue
inductively, following the described procedure. Note that
$\l^1= \l'_1 \op \a_1$ and $\m_+$ admit the following decompositions
as direct sum of linear spaces, each term of which is a Lie subalgebra:
\beqa
&&\l^1=\left(\l^1 \cap \Ad_{\dot{w}^1}^{-1}(\n^1_+) \right) \ds
       \left(\l^1 \cap \Ad_{\dot{w}^1}^{-1}(\g^1)   \right) \ds
       \left(\l^1 \cap \Ad_{\dot{w}^1}^{-1}(\n^1_-) \right),
\label{l^1decomp} \\ 
&&\m_+=\left(\ol{\n}^1_+ \cap \Ad_{\dot{w}^1}^{-1}(\n^1_+) \right) \ds
       \left(\m_+        \cap \Ad_{\dot{w}^1}^{-1}(\g^1)   \right) \ds
       \left(\ol{\n}^1_+ \cap \Ad_{\dot{w}^1}^{-1}(\n^1_-) \right).
\label{n^1+decomp}
\eeqa
These formulas imply that the product 
\beqa
&& \Big( L^1        \cap \Ad_{\dot{w}^1}^{-1}(N^1_+)  \Big)
   \Big( \ol{N}^1_+ \cap \Ad_{\dot{w}^1}^{-1}(N^1_+)  \Big)
   \Big( L^1        \cap \Ad_{\dot{w}^1}^{-1}(G^1)    \Big)
   \Big( M_+        \cap \Ad_{\dot{w}^1}^{-1}(G^1)    \Big)
\times
\nn \\
&& \times \; \Big( L^1  \cap \Ad_{\dot{w}^1}^{-1}(N^1_-)  \Big)
   \Big( \ol{N}^1_+     \cap \Ad_{\dot{w}^1}^{-1}(N^1_-)  \Big)
\label{product}
\eeqa
is a dense, open subset of $G_+,$ because the linear span of 
the Lie algebras of all factors is $\g_+.$ Although
$\ol{\n}^1_+ = \n^1_+(=\n_+)$ and 
$\ol{N}^1_+ = N^1_+(=N_+),$
we use both notations since this is how the groups
$\ol{N}^k_+ \neq N^k_+$ appear at the next steps $k \geq 2.$

  {} From this it follows that the double coset 
$[(l \dot{w}, e)]= [(g^1 \dot{w}^1, e)]$ of
$\D(G_+)$ has a dense, open subset consisting of elements
of the form
\beq
(m' l' g^1  \dot{w}^1 l''_+ n''_+ l''_\st m''_\st l''_- n''_-,
\pr_+(l') \pr_+(l''_+) \pr_+(l''_\st)  \pr_+(l''_-) ) 
\label{1coset}
\eeq
for some arbitrary elements $l''_+,$ $n''_+,$ $l''_\st,$ $m''_\st,$
$l''_-,$ and $n''_-$ of the groups of the product
\eqref{product} (taken in the same order). Note that after commuting
the terms $l''_+$ and $n''_+$ with $\dot{w}^1$ in the first component, they
will be absorbed by the term $m' \in M_+.$

The group $L^1$ by itself has a dense, open subset
\[
  \left(L^1   \cap \Ad_{\dot{w}^1}^{-1}(N^1_-)  \right)
  \left(L^1   \cap \Ad_{\dot{w}^1}^{-1}(G^1)    \right)
  \left(L^1   \cap \Ad_{\dot{w}^1}^{-1}(N^1_+)  \right)
\]
If the term $l'$ from \eqref{1coset} is taken in it, 
after some manipulations one shows that
the coset $[(g^1 \dot{w}^1, e)]$ has the following dense, open subset
\beqa 
&&       \Big( N^1_+, e \Big) \:
         \Big( (\id, \pr_+)(N^2_-)  \Big) \:
         \Big(e, \pr_+(N^2_+)  \Big) \times
\nn \\
&&\times \: \Big( H^\ort_1 N^2_+, e \Big) \:
            \Big( (\id, \pr_+)(L^2) \Big) \:
            \Big( g^1, e \Big) \:
            \Big( (\Ad_{\dot{w}^1}, \pr_+)(L^2) \Big) \:
            \Big( \ol{N}^2_+ \Ad_{\dot{w}^1}(H_1^\ort), e \Big)
\times
\nn \\
&&\times \: \Big( (\id, \pr_+ \Ad_{\dot{w}^1}^{-1})
                  (\Ad_{\dot{w}^1}(L^1) \cap N^1_-) \Big) \:
            \Big( \Ad_{\dot{w}^1}(\ol{N}^1_+) \cap N^1_-, e \Big) \:
            \Big( \dot{w}^1, e \Big).
\label{2coset}
\eeqa

Here for a subgroup $R_1$ of $R_2,$ and two endomorphisms $\sig_1,$
$\sig_2$ of $R_2,$ $(\sig_1, \sig_2)(R_1)$ denotes the image of the
composition of the diagonal embedding of $R_1$ in $R_2 \times R_2$
composed with $\sig_1 \times \sig_2.$

An important point
in the derivation of formula \eqref{2coset} is that $H_1^\ort$
commutes with $L^1 = L'_1 A_1$ (see \eqref{1.ortcompl})  and thus
with $N^2_- \sub L^1.$ Later this will be used for all groups 
$N^k_- \sub L^1,$ $k =2, \ldots, K.$

Note that $\La_1 \sub Z(L'_1 \, A_1) \sub H$ implies
$\La_1 \sub L^k$ $(\forall k).$ The product on the second line
of \eqref{2coset} is the double --
left $(H^\ort_1 N^2_+) \, L^2$  and right
$\Ad_{\dot{w}^1}(L^2) \, (\ol{N}^2_+ \Ad_{\dot{w}^1} (H_1^\ort))$
coset of $(g^1, e)$ in $G^1 \times (L^2/ \La_1)$ where
the two groups are embedded in
$G^1 \times (L^2/\La_1)$ as it is shown in the formula.
Using the representation \eqref{wk+1} of $g^1$ and \eqref{compute}
one easily shows that $[(g^1 \dot{w}^1, e)] = [(g^2 \dot{w}^2, e)]$
and procceeds by induction on $k.$

This procedure proves that $[(l \dot{w}, e)] = [(g^K \dot{v}, e)]$
(recall \eqref{v}). To prove the rest of \leref{lemma1}, let us
rearrange the coset from the last step similarly to \eqref{comp1}:
\beqa
&& \Big(H^\ort_1 N^K_+,e \Big) \/
   \Big( (\id, \pr_+)(L^K) \Big) \/
   \Big( g^K,e \Big) \/
   \Big( (\Ad_{\dot{v}}, \pr_+)(L^K) \Big) \/
   \Big( \ol{N}^K_+ \Ad_{\dot{v}}(H^\ort_1), e \Big) 
\label{cos-expl} \\
&&=\Big( N^K_+,e \Big) \/
   \Big( (\id, \pr_+) (L^K) \Big) \/
   \Big(\wt{TC}^{\dot{v}}_{L^v \times H^\ort_1 \times H^\ort_1}(g^K),
        e \Big)
\label{cos-last}
\eeqa
where we used that $\ol{N}^K_+ = N^K_+$ normalizes 
$\ol{G}^K =G^K.$
The action $\wt{TC}^{\dot{v}}$ of
$L^v \times H^\ort_1 \times H^\ort_1$ on $G^v$ is defined by
\beq
\wt{TC}^{\dot{v}}_{(l, h', h'')}(g) =
        h' l^{-1} g \Ad_{\dot{v}}(l h''), \;
        \mbox{for} \; g \in G^v, \ l \in L^v,
        \; h', \: h'' \in H^\ort_1.
\label{wt-tc}
\eeq
In \eqref{cos-last} 
the term $i_+(H^\ort_1 \, L^K)$
is replaced by $i_+(H^\ort_1) i_+(L^K)$ and
$(\Ad_{\dot{v}}\times \id)(i_-(L^K \, H^\ort_1))$ by
$(\Ad_{\dot{v}}\times \id)(i_-(L^K)i_-(H^\ort_1)).$
(The intersection $H^\ort_1 \cap L^K=H^\ort_1 \cap L'_1 A_1$ 
is nontrivial in general.) 

In the case when $H^\ort_{1,2}$ are trivial, the action
\eqref{wt-tc} gives a simpler classification of the $i_+(G_+)$
double cosets (see \coref{cosetfull}).

Using the decomposition \eqref{fc} of $g^K$ the middle set
in \eqref{cos-last} can be rewritten as
\beq
\wt{TC}^{\dot{v}}_{L^v \times H^\ort_1 \times H^\ort_1}(g^K) =
TC^{\dot{v}}_{(G^v)'}(f^K) \: (c^K
\: \Cong^v(L^v) H^\ort_1 \Ad_{\dot{v}}(H^\ort_1)).
\label{gtofc}
\eeq

The inductive argument shows that the coset $[(l \dot{w}, e)]$ has the
following dense, open subset:
\beqa
&&  \Big( N^1_+, e \Big) \/
    \Big( \prod_{k=2}^K E^k \Big) \/ 
    \Big( (\id, \pr_+) (L^K) \Big) \/ \times
\nn \\
&& \times \/ \Big( TC^{\dot{v}}_{(G^v)'}(g^K) \/ 
   (c^K \: \Cong^v(L^v) H^\ort_1 \Ad_{\dot{v}}(H^\ort_1)), e \Big) \/
   \times
\nn \\
&&\times \/
  \Big( \prod_{k=K-1}^1 F_k \Big) \/
  \Big(\dot{w}^K \ldots \dot{w}^1, e \Big).
\label{3coset}
\eeqa
where the subsets $E^k$ and $F^k$ of $G$ are defined by
\[
E_k=\Big( (\id, \pr_+)(N^k_-) \Big) \:
    \Big( N^k_+ \times \pr_+(N^k_+) \Big)
\]
and
\beqa
F_k&=&\Big( (\Ad_{(\dot{w}^K \ldots \dot{w}^k)},
             \Ad_{(\dot{w}^K \ldots \dot{w}^{k+1})} \pr_+)
      (L^k \cap \Ad_{\dot{w}^k}^{-1}(N^k_-)) \Big) \times
\nn \\
&&\times \:
      \Big( \Ad_{(w^K \ldots w^{k+1})}
      (\Ad_{\dot{w}^k}(\ol{N}^k_+) \cap N^k_-), e \Big).
\nn
\eeqa
The order of the factors in the first product in \eqref{3coset}
is $k=2, \ldots, K$ and in the second one is $k=K-1, \ldots, 1.$

Note that $\Ad_{\dot{w}^1}(\ol{\n}_+^1) \cap \n_-^1$ and
$\Ad_{\dot{w}^1}(\l^1) \cap \n_-^1$
are subalgebras of $\n_-^1$ and thus do not intersect 
$\n_+^1$ and $\g^1.$ Using this, the fact that $\pr_\pm$ are one to one 
on any unipotent subgroup of $L'_{1,2} \, A_{1,2},$ and the Bruhat Lemma 
we see that if two subsets of $G$ of the type \eqref{2coset} are equal,
then they are related to one and the same $w^1 \in W$ and
the corresponding double cosets of
$G^1 \times (L^2/\La_1)$ (on line two) coinside.
To finish the proof of \leref{lemma1}, assume
that two double cosets of the type $[(g^K \dot{v}, e)]$ obtained from
described procedure coincide.
Then their dense, open subsets from \eqref{2coset} will have nontrivial
intersection. Using the above argument repeatedly we get that
the two cosets are related to one and the same elements
$w^K,$ $\ldots,$ $w^1 \in W$ (and thus the same $v\in V)$ and
the corresponding to them
orbits $\wt{TC}^{\dot{v}}_{L^v \times H^\ort_1 \times H^\ort_1}(g^K)$
in $G^v$ coincide. The last is equivalent to saying that the orbits
$TC^{\dot{v}}_{(G^v)'}(f^K)$ corresponding to the two double cosets 
and the images of $c^K$ in the factor group \eqref{abelf}
are equal modulo the multiplication action of $\La^v_1$
on $(G^v_1)' \times Z(G^v)^\ci.$ \hfill \qed 
\\ \hfill \\
{\em{Proof of \leref{lemma2}.}} Let $v^0$ be an arbitrary
element of $V.$ Decompose it as
$v^0 = v^1 w^1 u^1$ with $v^1, u^1 \in W_1$ and 
$w^1$ the minimal length representative in $W$ of the 
double cosets $W_1 v W_1.$ Since $v^0$ is minimal
length representative of a coset from $W/W_1,$
using part 3 of \leref{MinLe} we get that
$u_1$ is the identity element of $W_1.$ 
Define two reductive Lie subalgebras $\g^2$ and $\ol{\g}^2$
of $\g^1$ by \eqref{gk+1}, \eqref{bargk} with $k=1$ and let
$W^2,$ $\ol{W}^2$ be the corresponding subgroups
of $W.$ The minimality property of $v^0$ implies
that for any positive root $\al^1$ of $\g^1,$
$v(\al^1)$ is a positive root of $\g.$ This means that
for any positive root $\al^2$ of $\ol{\g}^2,$
$v^1(\al^2)$ is a positive root of $\g^1.$ Therefore
$v^1$ is minimal length representative in $W^1$ of
a coset in $W^1 / \ol{W}^2.$ Continuing by induction we 
define $w^k,$ $v^k,$ $\g^k,$ $\ol{\g}^k,$ $W^k,$ $\ol{W}^k.$
At some step $k=K+1$ we will have $\g^{K+1} = \ol{\g}^{K+1}$
which is the same as
$\g^K = \Ad_{(\dot{w}^K \ldots \dot{w}^1)} (\g^K).$
As in the proof of \leref{lemma1} set
\[
v = w^K \ldots w^1
\]
and then
\beq
v^0 = v^K v.
\label{vK}
\eeq
We already proved that $\g^v= \g^K.$ From \eqref{vK} it follows that
$\g^{v^0} = \g^K,$ so both $v$ and $v^0$ have the property
that they are minimal length representatives of one and the 
same coset from $W^K \backslash W / W^1.$ The uniqueness of minimal
length representative implies $v^0 = v$
(see part 1 of \leref{MinLe}). Comparing this procedure
with the one from the beginning of the proof of \thref{G+coset}
completes the proof of
the existence part of \leref{lemma1}.
The uniqueness part is obtained by repatedly applying
part 1 of \leref{MinLe}.
\hfill
\qed 
\\ \hfill \\
{\em{Proof of formula \eqref{dim}.}}
The double coset $[(g^K \dot{v}, e)]$ can be considered as an orbit of
an action of the $2 \dim \g_+$ dimensional group $G_+ \times G_+$ on
$\D(G_+).$
	At the $k$-th step the set
\[
(\ol{N}_+^k \cap \Ad_{\dot{w}^k}^{-1}(N_+^k),e)
\]
coming from the second factor is annihilated $(k=1, \ldots, K-1).$
Its dimension is
\[
\dim \n^k_+ - \dim \n^{k+1}_+ - l(w^k).
\]
Using part 2 of \leref{MinLe} we obtain that the sum of the
dimensions of the annihilated sets is
\[
\dim \n_+ - \sum_{k=1}^K l(w^k) = \dim \n_+ - l(v).
\]
At the end one needs to subtract also the dimension
of the stabilizer of $g^K$ of the action
$\wt{TC}^{\dot{v}}$ of 
$L^v \times (H^\ort)^{\times 2}.$ 
Passing to dimensions of orbits and using 
formula \eqref{gtofc} we see that
\[
\dim \Stab_{\wt{TC}^{\dot{v}}} (g^K) =
\dim \l^v + 2 \dim \h^\ort_1 - TC^{\dot{v}}_{(G^v)'}(f^K)
- \dim( \Cong^v(L^v) \: H^\ort_1 \: \Ad_{\dot{v}}(H^\ort_1) ).
\]
Formula \eqref{dim} is obtained by combining
the above results with the fact that two elements
of the set \eqref{3coset} are equal if and only if they have
equal terms in each factor of the product.
\\ \hfill \qed 
\\ \hfill \\
{\em{Identification of the data of \thref{G+coset} corresponding
to different choices of a representative of $v \in V$ in $N(H)$:}}
If $\dot{v}$ and $\ddot{v}$ are two different representative of
$v \in V$ in $N(H)$ then
\[
\ddot{v} = h \dot{v}, \; \mbox{for some} \; h \in H,
\]
since the centralizer of a maximal torus in a complex reductive
group coincides with the torus. Let $h = h^v h^\ci$ for some
$h^v \in (G^v)'$ and $h^\ci \in Z(G^v)^\ci.$ A direct computation
shows that the right multiplication by $h^v$ in $(G^v)'$
maps $TC^{\ddot{v}}$ orbits to $TC^{\dot{v}}$ orbits: 
\[
TC^{\ddot{v}} (g h^v) = TC^{\dot{v}}(g) \, h^v, \;
\mbox{for all} \; g \in G^v.
\]
Since $\La_1^v$ is in the center of $(G^v)' \times Z(G^v)^\ci$ the right
multiplication of $h^v$ on the first component of $(G^v)' \times
Z(G^v)^\ci$ commutes with the regular action of $\La_1^v.$ This
identifies the data in \thref{G+coset} corresponding to two different
choices of representatives $\dot{v}$ and $\ddot{v}$ of $v$ in $N(H).$
\hfill \qed

\bco{cosetfull}If the Poisson structure on $G$ comes from an $r$-matrix
$r$ for which $\g_+ \sup \h,$
the $i_+(G_+)$ double cosets of $\D(G_+)$ are classified
by a minimal length representative $v$ of a coset in $W/W_1$
and an orbit of the twisted
conjugation action of $G^v$ on $G^v/\La_1$ defined by
\beq
\wt{TC}^{\dot{v}}_{g} (f) =
g^{-1} f \Ad_{\dot{v}}(g); \; g \in G^v, \; f \in G^v/ \La_1.
\label{wttc}
\eeq
The dimension of the coset corresponding to the $\wt{TC}^{\dot{v}}$ orbit
through $f \in G^v/\La_1$ is
\beq
\dim \g_+
+ \dim \l_1 - \dim \l^v 
+l(v) + \wt{TC}^{\dot{v}}_{G^v}(f).
\label{dimfull}
\eeq
\eco

The action \eqref{wttc} is the push down of the twisted conjugation
action \eqref{wt-tc} of $G^v$ on itself to the factor group $G^v/\La_1.$
It is correctly defined since $\La_1 \sub Z(G^v).$ 
\\ \hfill \\
\proof In the considered case the groups $H^\ort_{1,2}$ are trivial
and $L^v = G^v$ (see \reref{fullHgr}). 
The statement follows directly from the proof of \thref{G+coset}, in
particular \eqref{cos-last}.
\hfill \qed 

\bex{Cosetsid} The simplest coset in $W/ W_1$ is 
$W_1,$ whose minimal length representative in $W$ is 
the identity element $\id$ of $W.$ 
In this example we discuss the $i_+(G_+)$ double cosets
corresponding to it. We will assume that $\id$ is represented in 
$N(H)$ by the identity element $e$ of $G.$ 
The group $L^\id$ is simply $L'_1 \, A_1$ and
for any coset of this type
the inductive procedure of the proof of \thref{G+coset}
has only one step, i.e. $K=1.$ In the product \eqref{product} 
all groups except the second, the third, and the forth 
are trivial. Besides this the full group $G_+$ is represented as
\[
G_+ = N_+ \/ (L'_1 \: A_1) \/ H^\ort_1,
\]
not just a dense subset of it. (This is nothing but \eqref{Gpm}.)
The expression \eqref{3coset} reduces to:
\beq
[(f c, e)]
= \Big( N_+, e \Big) \:
\Big( (\id, \pr_+)(L'_1 A_1) \Big) \:
\Big( TC^e_{L'_1}(f) \: (c \, H^\ort_1), e \Big)
\label{idfc}
\eeq
for all $f \in L'_1$ and $c \in Z(L_1)^\ci.$
The action $TC^e $ is the usual conjugation action of $L'_1$ on itself
\[
TC^e_{l} (f) = l^{-1} f l, \; l, f \in L'_1.
\]
Note also that $\Cong^\id(L^\id) H^\ort_1 \Ad_e(H^\ort_1) = H^\ort_1.$

The subgroup $\La^\id_1$ of $L'_1 \times Z(L_1)^\ci$
is the pull back of $\La_1$ under the
map $L'_1 \times Z(L_1)^\ci \ra L'_1 Z(L_1)^\ci = L_1.$

Thus the cosets corresponding to $\id \in V$ are classified by
the conjugation orbits on $L'_1$ and the elements of the abelian
group $Z(L_1)^\ci / H^\ort$ taken modulo the multiplication action
of $\La^\id_1$ on $L'_1 \times Z(L_1)^\ci.$

The dimension of the coset \eqref{idfc} is
\[
\dim(\g_+) +  
   \dim(TC^e_{L'_1} (f))
\]
which follows from \eqref{dim} after some cancellations.
\eex
\bex{cosetstan}Consider the case of the standard Poisson structure on
a semisimple Lie group $G$ (see 
Examples \ref{e1.stand} and \ref{e2.stand}).

The groups $W_{1,2}$ are trivial, so $W/ W_1= W.$ Any $i_+(B_+)$ double
coset of $D(B_+) \cong G \times H$ is equal to a coset of the type 
$[(h \dot{w}, e)],$ with $w \in W,$ $h \in H.$ For all $w \in W$
the procedure from the proof of \thref{G+coset} again has only one step
and $L^w =H.$ All groups in the product \eqref{product} except 
the second, the third, and the sixth are trivial and it reduces to 
\[
B_+ = (N_+ \cap \Ad^{-1}_{\dot{w}}(N_+)) \:
      H \: (N_+ \cap \Ad^{-1}_{\dot{w}}(N_-)).
\]
Since the full group $G_+=B_+$ is represented in this way,
using \eqref{3coset} we can get an explicit expression for
any double coset of $\D(B_+):$ 
\beq
[(h \dot{w}, e)] =
(N_+, e) \: d(H) \: ( h \: \Cong^w(H), e ) \:
(\Ad_{\dot{w}}(N_+) \cap N_-, e) \:
(\dot{w}, e), \: \forall
w\in W, \: h \in H
\label{stand-cos}
\eeq
where
\[
\Cong^w(H) = \{(h')^{-1} \Ad_{\dot{w}}(h') | \: h' \in H \}
\]
and $d(H)$ is the image of the diagonal embedding of $H$ in
$G \times H.$
The set $h \Cong^w(H)$ is the $\wt{TC}^{\dot{w}}_H$ orbit
through $h$ from \coref{cosetfull}. The general twisted conjugation 
orbits in
\thref{G+coset} and \coref{cosetfull} can be considered as nonabelian
analogs of this action. 

The different cosets \eqref{stand-cos} are classified by the elements
$w$ of the Weyl group and of the abelian group $H / \Cong^w(H)$ (or
equivalently the orbits of the $\wt{TC}^{\dot{w}}$ action of $H$ on
itself). According to \thref{int} this data classifies
the symplectic leaves of $i_-(B_-)/\Sig.$

Both \eqref{stand-cos} and \eqref{dimfull} give 
\[
\dim [(h\dot{w}, e)] = \dim \b_+ + l(w) + \dim (\Cong^w(H)).
\]

The above description of $i(B_+)$ double cosets of $\D(B_+)$  
(and symplectic leaves of $i_-(B_-)/\Sig_-)$ for
the standard Poisson--Lie structure on $G$ agrees with the
results in \cite{DKP, HKKR}.
\eex
\bex{CGcosets} The Cremmer--Gervais Poisson structure on $SL(n+1)$
was reviewed in Examples \ref{eCG} and \ref{eCGgroup}. In this
case $W \cong S^{n+1}$ and $W_i$ can be identified with
the subgroups
$S_i^{n+1}$ of $S^{n+1}$ 
that permute the first and the last $n$ letters.
Denote by $H$ the maximal torus of $SL(n+1)$ consisting
of diagonal matrices.
We will assume that any element of $S^{n+1}$ is represented
in the normalizer of $H$ by a multiple of the corresponding
standard permutation matrix. Any minimal length representative $\sig$
in $S^{n+1}$ of a coset from $S^{n+1} / S^n_1$ should have the
property that $\sig(1),$ $\ldots,$ $\sig(n)$ is an increasing sequence
(see part 1 of \leref{MinLe}). From this one finds that the set
of such representatives consists of
$\sig^j_1 \in S^{n+1},$ $j=0, \ldots, n,$ defined by
\[
\sig^j_1(p) = p , \, p=1, \ldots, j;  \;
\sig^j_1(p) = p + 1, \, p=j+1, \ldots, n; \; 
\sig^j_1(n+1) = j+1.
\]
In a compact form, they are the permutations
\beq
\sig^j_1 = ( (j+1) \, (j+2) \, \ldots \, (n+1)). 
\label{sig1}
\eeq
The groups $G^{\dot{\sig}^j_1}$ are
\[
G^{\dot{\sig}^j_1} = \{ \diag(A, a_{j+1}, \ldots, a_{n+1}) |\,
A \in GL(j), a_p \in \Cset; \;
\det(A) a_{j+1} \ldots a_{n+1} =1 \} \sub SL(n+1).
\]
The twisted conjugation action of $G^{\sig^j_1}$ on itself
is given by
\beqa
&&\wt{TC}^{\dot{\sig}^j_1}_{\diag(A, a_{j+1}, \ldots, a_{n+1})} \,
\diag(B, b_{j+1}, \ldots, b_{n+1}) =
\nn \\
&&\diag(A^{-1} B A,
a_{j+1}^{-1} b_{j+1} a_{\mu(j+1)},
\ldots, a_{n+1}^{-1} b_{n+1} a_{\mu(n+1)})
\eeqa
where $\mu = (\sig^j_1)^{-1}.$
To each such orbit one can associate the conjugation
orbit of $GL(j)$ trough $B.$ The correspondence
is bijective and the difference in their dimensions is
$n-j.$ Note also that $l(\sig^j_1)= n-j.$

Thus in the case of the Cremmer--Gervais structure
the $i_+(G_+)$ double cosets of $\D(G_+),$
or in other words the symplectic leaves of
$i_-(G_-)/\Sig_-$ (recall \thref{int}) are classified
by a choice of a number $j=0, \ldots, n$ and a conjugation
orbit on $GL(j).$ The dimension of the corresponding
double coset
is
\[
n(n+1) + (n-j)(n+j+1) + \dim \{ \Ad_A(B) | \,
A \in GL(j) \}.
\]
\eex
%%%%%%%%%%%%%%%%%%%%%%%%%%%%%%%%%%%%%%%%%%%%%%%%%%%%%%%%%%%%%%%%%%%%%%%%%%%%%%%
\subsection{Structure of the group
$\Sig_- = i_+(G_+) \cap i_-(G_-) \sub \D(G_+)$}
Denote
\[
\Sig_- = i_+(G_+) \cap i_-(G_-) \sub \D(G_+).
\]
It is a discrete subgroup of $\D(G_+)$ since
$\Lie(i_+(G_+)) \cap \Lie(i_-(G_-)) = 0.$ 

\bth{Sig+}The group $\Sig_-$ is the intersection of the $i_\pm$ images
of the tori $(L'_i\cap H) A_i \sub G_\pm$ in the maximal torus
$H \times ((L'_1\cap H) A_1/ \La_1)$ of the reductive group
$G \times (L'_1  A_1/\La_1) \cong \D(G_+).$
\eth
\proof Using the recursive procedure from Subsect.~3.3 we will show that
the statement holds for $G^{(k)},$ assuming its validity for $G^{(k+1)}.$
Note that it is trivial for $G^{(\ord \tau)}$ since $G^{(\ord \tau)} = H.$ 
Clearly we can restrict ourselves to the step from $G^{(1)}$ to
$G=G^{(0)}.$ 

Let
\[
l_i \in L'_i A_i, \, i=1,2;  \; m_\pm \in M_\pm \; 
\]
be such that
\beq
(l_1 m_+, \pr_+(l_1))=(m_- l_2, \Th^{-1} \ci \pr_-(l_2)) 
\in i_+(G_+) \cap i_-(G_-).
\label{inter1}
\eeq
   {} From the first components we get
\beq
l_1 m_+ = m_- l_2 \in
\left((L'_1 \: A_1) \: M_+ \right) \cap
\left((L'_2 \: A_2) \: M_- \right).
\label{inter2}
\eeq
Let $N^{(1)}_{1\pm}$ denote the unipotent subgroups
of $L'_1$ generated by one parameter subgroups
of positive / negative roots of $\l'_1$ that are not roots
of $\l'_2.$ Similarly we define $N^{(1)}_{2 \pm}$ 
as the unipotent subgroups of $L'_2$ corresponding to roots
of $\l'_2$ which are not roots of $\l'_1.$
Using the standard
fact that the intersection of the two parabolic subgroups
$P_+$ and $P_-$ of $G$ is connected we get
\beqa
&&m_+ \in H_1^\ort \ltimes N^{(1)}_{2+},
\label{prodm+} \\
&&m_- \in N^{(1)}_{1-} \rtimes H_2^\ort,
\label{prodm-}
\eeqa
and
\beqa
&&l_1 \in N^{(1)}_{1-} \rtimes 
\left( (L'_1 \cap L_2) \: A_1 \right), 
\label{prodl1} \\
&&l_2 \in \left( (L'_2 \cap L_1) \: A_2 \right) 
\ltimes N^{(1)}_{2+}. 
\label{prodl2}
\eeqa

Let $m_+ = h_1 n_{2+}$ and $l_1 = n_{1-} l^{(1)}_1$ be the decompositions
of $m_+$ and $l_1$ according to the semidirect products \eqref{prodm+}
and \eqref{prodl1}. Eq. \eqref{inter2} implies that the decompositions of
$m_-$ and $l_2$ according to \eqref{prodm-} and \eqref{prodl2} should be
$m_-= n_{1-} h_2,$ $l_2= k^{(1)}_2 n_{2+}$ for some $h_2 \in H^\ort_2$ and
$k^{(1)}_2 \in (L'_2 \cap L_1) \: A_2$
such that
\beq
l^{(1)}_1 h_1 = h_2 k^{(1)}_2. 
\label{l1k2}
\eeq
Observe that 
\beqa
&&(L'_1 \cap L_2) \: A_1 = L^{(1)}_1, \; N^{(1)}_{2+} = N^{(1)}_+,
\label{groups} \\
&&\th(\Lie((L'_1 \cap L_2) \, A_1))= \Lie(L^{(1)}_2 \, A^{(1)}_2),
\; \th(\Lie(N^{(1)}_{1-})) = \Lie(N^{(1)}_-).
\label{algth}
\eeqa
Since $N^{(1)}_-$ is a unipotent subgroup of $L'_2 \, A_2$
and the kernel $\La_2$ of the projection \\
$\pr_-: \, L'_2 \, A_2 \ra L'_2 \, A_2 /\La_2$
consists of semisimple elements, eqs. \eqref{algth} 
imply the existence of elements $n^{(1)}_- \in N^{(1)}_-$ and 
$l^{(1)}_2 \in (L^{(1)})' \, A^{(1)}$ such that
\beqa
&&\Th \ci \pr_+(n_{1-}) =\pr_-(n^{(1)}_-),
\label{le1} \\
&&\Th \ci \pr_+(l^{(1)}_1) =
        \pr_-(l^{(1)}_2).
\label{le2}
\eeqa 
Comparing the second components of \eqref{inter1} we get
$\Th \ci \pr_+(n_{1-} l^{(1)}_1) = \pr_-(k^{(1)}_2 n_{2+}).$
Combining this equation with \eqref{le1} and \eqref{le2} gives
\[
\pr_-( n^{(1)}_- l^{(1)}_2)= \pr_-(k^{(1)}_2 n_{2+}).
\]
By modifying $l^{(1)}_2 \in (L^{(1)})' \, A^{(1)}$
with an element of $\La_2 \sub (L^{(1)})' \, A^{(1)}$
one can make it satisfy
\beq
n^{(1)}_- l^{(1)}_2 = k^{(1)}_2 n_{2+}
\label{eq1}
\eeq
while keeping \eqref{le2}. Expressing $k^{(1)}_2$ in terms
$l^{(1)}_1$ from \eqref{l1k2} and substituting it in
\eqref{eq1} gives
\beq
h_2 n^{(1)}_- l^{(1)}_2 = l^{(1)}_1 h_1 n_{2+}.
\label{eq2}  
\eeq 

Because $l^{(1)}_1 \in (L'_1 \cap L_2) \: A_2  = L^{(1)}_1$
and $n_{2+} \in N^{(1)}_{2+} = N^{(1)}_+$ eqs. \eqref{eq2} and
\eqref{le2} imply that
\[
(l^{(1)}_1 h_1 n_{2+}, \pr_+(l^{(1)}_1) ) =
(h_2 n^{(1)}_- l^{(1)}_2, \Th^{-1} \ci\, \pr_-(l^{(1)}_2)) 
\in i^{(1)}_+(G^{(1)}_+) \cap i^{(1)}_-(G^{(1)}_-) \sub \D(G^{(1)}_+).
\]

Using the inductive assumption for validity of the statement
of the Lemma for $G^{(1)},$ we get that $l^{(1)}_1 \in H,$ 
$n_{2+} = e$ and therefore
$(l_1 m_+, l_1) \in H \times ((L'_1 \cap H) \times A_1).$
\hfill \qed 
\\ \hfill 

As a consequence of \thref{Sig+} we also obtain that the $\Sig_-$--groups
associated to the Poisson--Lie groups $G,$ $G^{(1)},$ $\ldots,$
$G^{(\ord \tau)}$ are naturally isomorphic.

   {} From \thref{Sig+} one easily derives 
an explicit formula for $\Sig_-$ in
terms of the lattices of kernels of exponential maps for certain tori.
We will restrict ourselves to the case of an $r$-matrix
for which $\h^\ort_i$ are trivial since in the general case
the formula is too cumbersome.

Denote by $\ker$ and $\ker'$ the kernels of the exponential maps
for the tori $H$ and $H/ \La_2.$

\bco{Sig}In the case when the $r$-matrix $r$ is such that
$\g_+ \sup \h,$ the group $\Sig_-$ is isomorphic to
\beq
\ker' / (\ker' \cap (1-\th) \ker). 
\label{Sigma}
\eeq
If in addition the groups $\La_i$ are trivial, then 
\beq
\Sig_- = \{ (h, h) | \, h \in H, \, \Th(h) = h \}.
\label{Sinv}
\eeq
\eco

One easily deduces from the assumption $\g_+ \sup \h$
and the fact $\D(\g_+) = i_+(\g_+) \ds i_-(\g_-),$ 
that $1-\th$ is a linear automorphism of $\h.$
This implies that $(1-\th)^{-1}\ker'$ is a lattice in $\h.$
The quotient \eqref{Sigma} of the two lattices 
is not a finite group in general.
\\ \hfill \\ 
{\em{Proof}} of \coref{Sig}. According to \thref{Sig+} $\Sig_-$ is the
subgroup of $H \times (H/ \La_1)$ consisting of elements
\beq
(l_1, \pr_+(l_1)) = (l_2, \Th^{-1} \ci \pr_-(l_2) )
\label{Sigeq}
\eeq
for some $l_i \in H.$ Let $l_1 = exp(x)$ for some $x \in \h.$
The second component of eq. \eqref{Sigeq} gives
$l_2 = \exp( \th x) \la_2$ for some $\la_2 \in \La_2.$ From the first
component we obtain 
\[
\exp( (1 - \th) x) = \la_2 \in \La_2
\]
and thus $x \in (1-\th)^{-1} \ker'.$ Taking the quotient
\beq
(1-\th)^{-1} \ker' / ( (1-\th)^{-1} \ker' \cap \ker)
\label{1-th}
\eeq
cancels the elements $x \in \h$ for which $\exp(x)=e.$
The final answer is obtained after applying $(1-\th)$ to
\eqref{1-th}.

Eq. \eqref{Sinv} follows directly from \eqref{Sigeq}.
\hfill \qed

\bex{Sigstand} Consider the case of the standard Poisson
structure on a complex semisimple Lie group $G.$ Then $\La_2$ is trivial
and 
$\ker' = \ker.$ \coref{Sig} combined with the facts that
the Cayley transform $\th: \; \h \ra \h$ is $- \id|_\h$ and
$\La_2$ is trivial gives
\beq
\Sig_- \cong  \ker / (2 \ker) \cong
\Zset^{\rank \g} / (2 \Zset^{\rank \g})
\cong \Zset_2^{\rank \g},
\label{SigmaSt1}
\eeq
or equivalently
\beq
\Sig_- = \{ (h, h) | \, h \in H, \, h^2 = e \}.
\label{SigmaSt2}
\eeq
This coincides with the results in \cite{DKP, HKKR}.
\eex

\bex{SigCG} For the Cremmer--Gervais structure on $SL(n+1)$
the groups $\La_i$ are trivial. The group $\Sig_-$ is
explicitly given by
\beqa
\Sig_- &=& \{ (h, h) | \, h= \diag(h_1, \ldots, h_{n+1}); \;
h_{i+1} = h_i, \,  i=1, \ldots, n \}
\nn \\
&=&\{ (\e E, \e E) | \,\e \in \Cset, \, \e^{n+1} = 1 \} \cong \Zset_{n+1}
\label{SigmaCG}
\eeqa
where $E$ denotes the identity matrix of size $n+1.$
This can also be proved directly from \eqref{Sigma} using
the facts that the kernel of the exponential map on $\h$ is
\[
\ker = 2 \pi i \Zset \al_1 \op \cdots \op 2 \pi i \Zset \al_n.
\]
and the Cayley transform on Lie algebra level acts by
$\th(\al_j) = \al_{j+1}$ for $j=1, \ldots, n-1$
and $\th(\al_n)= -\al_1-\cdots - \al_n.$
\eex
%%%%%%%%%%%%%%%%%%%%%%%%%%%%%%%%%%%%
\subsection{Symplectic leaves of $G_-$}
In this Subsection we will summarize the results for the symplectic
leaves of $G_-$ and $i_-(G_-)/ \Sig_-$ obtained by 
combining the results from Subsections~4.1, 4.2 and Theorems \ref{t1.STS}, 
\ref{tint}. 

\bth{SymlG-} 1. The symplectic leaves of $i_-(G_-)/\Sig_-$ are
classified by the data from \thref{G+coset}: 
An element $v \in V,$ an orbit $TC^{\dot{v}}(f)$ in $(G^v)',$ and
an element of the abelian factor group 
$Z(G^v)^\ci / (\Cong^v(L^v) H^\ort_1 \Ad_{\dot{v}}(H^\ort_1)),$
all taken modulo the multiplication action of $\La^v_1.$
The dimension of the leaf corresponding to this data is:
\beq
\dim(\l'_1\op\a_1) - \dim \l^v
- \dim \h^\ort_1 +
l(v) + \dim TC^{\dot{v}}_{(G^v)'}(f)
+ \dim \Cong^v(L^v) H^\ort_1 \Ad_{\dot{v}}(H^\ort_1).
\label{leafdim2}
\eeq 
 
2. The leaves of $G_-$ are coverings of the leaves of $i_-(G_-)/\Sig_-$
under the projection map $G_- \ra i_-(G_-)/\Sig_-.$ The discrete group
$\Sig_-$ is characterized in \thref{Sig+}.
\eth

In a special case \coref{cosetfull} gives a simpler description of the
symplectic leaves of $G_-.$

\bco{symlG+full} In the case when $\g_+ \sup \h$ the symplectic
leaves of $i_-(G_-)/ \Sig_-$ 
are in one to one correspondence with the elements $v \in V$
and the orbits of the action $\wt{TC}^{\dot{v}}$ of $G^v$ on
$G^v/\La_1.$ {\em{(}}see \eqref{wttc}{\em{)}}. The dimension
of the leaf
corresponding to the orbit $\wt{TC}^{\dot{v}}(f)$ is
\[
\dim \l_1 - \dim \g^v + l(v)
+ \dim \wt{TC}^{\dot{v}}_{G^v}(f).
\]
The group $\Sig_-$ is explicitly computed in \coref{Sig}.
\eco

Note also 
that the results from Sect.~4.1--4.2 give some
information for the structure of the symplectic leaves 
of $i_-(G_-) / \Sig_-.$ In particular, the leaf
corresponding to the choice of the data specified
there  has a dense, open subset which is the intersection
of the set \eqref{3coset} (for $f^K =f,$ $c^K=c)$
with $i_-(G_-) i_+(G_+)$ projected onto $i_-(G_-) i_+(G_+)/ i_+(G_+).$
%%%%%%%%%%%%%%%%%%%%%%%%%%%%%%%%%%%%%%%%%%%%%%%%%%%%%%%%%%%%%%%%%%%%%%%%%%%%%%%
\sectionnew{Symplectic leaves of $G$}
This Section describes the set of symplectic leaves of the full group
$G.$ The proofs of all results will only be sketched since they go along
the lines of the proofs of  the previous Section.
We will restrict ourselves to the case when the derived subgroup
$G'$ of $G$ is simply connected. This is not an essential restriction
since any factorizable reductive Poisson--Lie group $R$ 
has a covering which is a Poisson--Lie group $G$ of the considered type. 
The leaves of $R$ are simply projections of the ones of $G.$

The above assumption implies in particular that
$L'_i$ are also simply connected and the restriction of the map
$\th$ to $\l'_1$ can be lifted to an isomorphism
\beq
\Th': \: L'_1 \ra L'_2.
\label{th'}
\eeq
It clearly satisfies
\[
\Th \ci \pr_+(l) = \pr_- \ci \Th' (l), \: \forall
l \in L'_1.
\]
%%%%%%%%
\subsection{$G^r$ double cosets of $\D(G)$}
%%%%%%%%
The double coset of an element $(g_1, g_2) \in G \times G \cong \D(G)$
will be denoted by
$[(g_1, g_2)]:$
\[
[(g_1, g_2)]:= G^r (g_1, g_2) G^r \sub G \times G.
\]
Similarly to \eqref{compute} one has
\beq
[(g_1 l, g_2)] =  [(g_1 , g_2 \: \Th' (l^{-1}) )], 
\;
[(l g_1 , g_2)] =  [(g_1 , \Th'(l^{-1}) \: g_2)],  
\forall l \in L'_1.
\label{compute12}
\eeq
As in Subsect.~4.1.1 it can be shown that any $G^r$ double coset
of $G \times G$ is equal to a coset of the form
\[
[(l_1 \dot{w}_1, \dot{w}_2 l_2)]
\]
where $l_i \in L_i$ and $w_i$ are representatives
in $W$ of cosets from $W_i \backslash W / W_i,$ $i=1, 2.$ 
For simplicity of notation it will be more convenient
to have the cosets written as
\[
[(l_1 \dot{w}_1, \dot{w}_2^{-1} l_2)]
\]
where $l_i$ and $w_i$ belong to the same sets as above.
%%%%%%
\subsubsection{Description of the set of double cosets}
The map $\tau$ from the Belavin--Drinfeld triple
$(\Ga_1, \Ga_2, \tau)$ associated to the $r$-matrix $r$ will
be assumed to be extended to a linear bijection between the set of
roots of $\l_1$ and $\l_2.$

Define the set $V$ as the subset of $W \times W$ consisting of pairs
of elements $(v_1, v_2)$ of $W$ which are representatives of minimal
length in $W$ of cosets from $W / W_i,$ $i = 1,2.$
For each such pair $(v_1, v_2)$ let $\g^{v_1, v_2}_1$ be the 
subalgebra of $\l_1$ containing $\h$ whose set of roots is stable
under
$v_1 \tau^{-1} v_2 \tau:$
\beqa
\g^{v_1, v_2}_1 = \span \{ h, \; \g^\al \; | \;
&&\mbox{for the roots} \; \al \; \mbox{of} \; \l_1 \; \mbox{such that}
\; (v_1 \tau^{-1} v_2 \tau)^n (\al) \; \mbox{is well defined}
\nn \\
&& \mbox{and is a root of} \; \l_1 \; \forall \; \mbox{integer} \; n \}.
\label{gwv1}
\eeqa
Note that $v_1 \tau^{-1} v_2 \tau$ is not defined 
on all root spaces $\g^\al$ of $\l_1$ but only on those 
for which $v_2 \tau (\al)$ is a root of $\l_2.$

In this notation the maximal subalgebra of $\l_2$ containing $\h$
whose set of roots is preserved by $v_2 \tau v_1 \tau^{-1}$ is
\beq
\g^{v_1, v_2}_2 = \Ad_{\dot{v}_2} \th (\g^{v_1, v_2}_1)' + \h.
\label{gw2}
\eeq

The minimality property of $v_i$ implies that 
$\g^{v_1, v_2}_i$ are reductive Lie subalgebras of $\l_i$
which are generated by $\h$ and root spaces of simple
roots of $\g.$ Denote by $G^{v_1, v_2}_i$ the connected
subgroups
of $G$ with Lie algebras $\g^{v_1, v_2}_i$ and by $W^{v_1, v_2}_i$ 
their Weyl groups, $i=1,2.$ Similarly to Sect.~4.1.2 one shows that
$v_i$ are minimal length representatives for the double
cosets $W^{v_1, v_2}_i \: v_i \: W_i.$
Clearly $\Ad_{\dot{v}_1} (\Th')^{-1} \Ad_{\dot{v}_2} \Th'$
is an automorphism of $(G^{v_1, v_2}_1)'.$

The connected subgroups of $G^{v_1, v_2}_i$ with Lie algebras 
$\g^{v_1, v_2}_i \cap (\l'_i \op \a_i)$ and \\  
$\g^{v_1, v_2}_i \cap \Ad_{\dot{v}_i} (\l'_i \op \a_i)$ 
will be denoted by $L^{v_1,v_2}_i$ and $\ol{L}^{v_1, v_2}_i.$ 
Note that
\beqa
&&\Ad_{\dot{v}_2} \th: \, \Lie(Z(L^{v_1, v_2}_1)) \ra
\Lie(Z(\ol{L}^{v_1, v_2}_2)) \; \mbox{and}
\nn \\
&&\th \Ad_{\dot{v}_1}^{-1}: \, \Lie(Z(\ol{L}^{v_1, v_2}_1)) \ra
\Lie(Z(L^{v_1, v_2}_2))
\nn
\eeqa
are linear isomorphisms.
For any $(v_1, v_2) \in V$ we will denote by $Z^{v_1, v_2}_i$
the following subgroups of the abelian group
$Z(G^{v_1, v_2}_1)^\ci \times Z(G^{v_1, v_2}_2)^\ci:$
\beqa
&& Z^{v_1, v_2}_1 = \{ (x, y) \in
Z(L^{v_1, v_2}_1)^\ci \times Z(\ol{L}^{v_1, v_2}_2)^\ci \: |
\Th \pr_+ (x) = \pr_- \Ad_{\dot{v}_2}^{-1} (y) \}^\ci,
\nn \\
&& Z^{v_1, v_2}_2 = \{ (x, y) \in
Z(\ol{L}^{v_1, v_2}_1)^\ci \times Z(L^{v_1, v_2}_2)^\ci \: |
\pr_+ \Ad_{\dot{v}_1}^{-1}(x) = \Th^{-1} \pr_-(y) \}^\ci.
\nn
\eeqa
Define also the following subgroup of
$Z(G^{v_1, v_2}_1)^\ci \times Z(G^{v_1, v_2}_2)^\ci$
\beq
Z^{v_1, v_2} = Z^{v_1, v_2}_1 \: Z^{v_1, v_2}_2 \:
(H^\ort_1 \Ad_{\dot{v}_1}(H^\ort_1) \times
H^\ort_2 \Ad_{\dot{v}_2}(H^\ort_2) ).
\label{Cv}
\eeq

\bth{Grcoset}The $G^r$ double cosets of $G \times G$ are classified 
by pairs $(v_1, v_2)$ of minimal length representatives in
$W$ of cosets from $W/ W_i,$
an orbit of the twisted conjugation action 
$TC^{\dot{v}_1, \dot{v}_2}$ of $(G^{v_1, v_2}_1)'$ on itself defined
by
\beq
TC^{\dot{v}_1, \dot{v}_2}_{g}(f) =
g^{-1} f \Ad_{\dot{v}_1} (\Th')^{-1} \Ad_{\dot{v}_2} \Th' (g);
\: g, \, f \in (G^{v_1, v_2}_1)', 
\label{tc2}
\eeq
and an element of the abelian factor group
\beq
(Z(G^{v_1, v_2}_1)^\ci \times Z(G^{v_1, v_2}_2)^\ci ) / 
Z^{v_1, v_2}.
\label{factab}
\eeq
This data has to be taken modulo the action
of the finite group
\beq
((G^{v_1, v_2}_1)' \cap Z(G^{v_1, v_2}_1)^\ci) \times
((G^{v_1, v_2}_1)' \cap Z(G^{v_1, v_2}_2)^\ci)
\label{cap}
\eeq
on $(G^{v_1, v_2}_1)' \times Z(G^{v_1, v_2}_1)^\ci
\times Z(G^{v_1, v_2}_2)^\ci$ where an element $(q_1, q_2)$
of \eqref{cap} acts on the direct product by multiplication by
$(q_1^{-1} \, \Ad_{\dot{v}_1} (\Th')^{-1}(q_2), q_1, q_2).$ 
{\em{(}}The choice of representatives of $v_i$ in the normalizer 
of the maximal torus $H$ of $G$ is not essential{\em{.)}}
\eth
{\em{An inductive procedure:}}
Let $[(l_1 \dot{w}_1, \dot{w}_2^{-1} l_2)]$ be an arbitrary $G^r$
double coset of $G \times G,$ where $l_i \in L_i$ and $w_i$ are
minimal length representatives of cosets from
$W_i \backslash W / W_i,$ $i=1, 2.$

We define two sequences of reductive Lie subalgebras $\g^k_i$
of $\g,$ elements $g^k_i \in \g^k_i,$ and $w^k_i \in W$ 
$(k \geq 1, \; i=1,2).$ Set $\g^1_i = \l_i,$ 
$g^1_i=l_i,$ $w^1_i=w_i.$ For convenience we also define
$\g^0_i = \g$ and {\em{extend $\th$ to a linear isomorphism
between $\l_1$ and $\l_2.$}}

As in the proof of \thref{G+coset} the defined objects will have
the following two properties:

(1) $\g^k_1$ and
\[
\ol{\g}^k_1:= \Ad_{(\dot{w}^{k-1}_1 \ldots \dot{w}^1_1)}
\th^{-1} (\g^k_2)
\]
are
reductive subalgebras of $\g^{k-1}_1 \sub \g$ generated by the 
Cartan subalgebra $\h$ and by root spaces of simple roots of $\g.$ 
The same property have the subalgebras $\g^k_2$ and
\[
\ol{\g}^k_2=\Ad_{(\dot{w}^{k-1}_2 \ldots \dot{w}^1_2)} \th (\g^k_1)
\]
of $\g^{k-1}_2.$

Let $G^k_i,$ $\ol{G}^k_i$ be the connected subgroups
of $G$ with Lie algebras $\g^k_i,$ $\ol{\g}^k_i$ and 
$W^k_i,$ $\ol{W}^k_i$ be their Weyl groups.

(2) $w^k_i$ are representatives of minimal length in $W^{k-1}_i$
of cosets from
$W^k_i \backslash W^{k-1}_i / \ol{W}^k_i.$ \\
If $\g^{k+1}_1 \neq \g^k_1$ or $\g^{k+1}_2 \neq \g^k_2$
we define $\g^{k+1}_1$ and $\g^{k+1}_2$ by
\beqa
&&\g^{k+1}_1 :=
\g^k_1 \cap \th^{-1} \Ad_{(\dot{w}^k_2 \ldots \dot{w}^1_2)}^{-1}(\g^k_2)
= \th^{-1} \Ad_{(\dot{w}^{k-1}_2 \ldots \dot{w}^1_2)}^{-1}
(\ol{\g}^k_2 \cap \Ad_{\dot{w}_2^k}^{-1}(\g^k_2) )
\label{g1k+1}, \\
&&\g^{k+1}_2 :=
\g^k_2 \cap \th \Ad_{(\dot{w}^k_1 \ldots \dot{w}^1_1)}^{-1} (\g^k_1)
= \th \Ad_{(\dot{w}^{k-1}_1 \ldots \dot{w}^1_1)}^{-1}
(\ol{\g}^k_1 \cap \Ad_{\dot{w}_1^k}^{-1}(\g^k_1) ).
\label{g2k+1} 
\eeqa

Denote with $N^{k+1}_{i \pm}$ and $\ol{N}^{k+1}_{i \pm}$ the unipotent
subgroups of $G^k_i$ generated by one parameter subgroups of 
positive / negative roots of $\g^k_i$ that are not roots of
$\g^{k+1}_i$ and $\ol{\g}^{k+1}_i.$ Set $\e_k = (-1)^k.$ Let us 
apply  the Bruhat Lemma to $g^k_1 \in G^k_1,$ $g^k_2 \in G^k_2$ and the
following pairs of parabolic subgroups:
$N^{k+1}_{1 \e_{k+1}} \rtimes G^{k+1}_1,$ 
$G^{k+1}_1 \ltimes \ol{N}^{k+1}_{1 +}$ 
of $G^k_1$ and 
$\ol{N}^{k+1}_{2 -} \rtimes G^{k+1}_2,$
$G^{k+1}_2 \ltimes N^{k+1}_{2 \e_k}$ of $G^k_2.$
We obtain
\beqa
&&g^k_1 = n'_1 g'_1 \dot{w}^{k+1}_1 g''_1 n''_1,
\label{w1} \\
&&g^k_2 = n'_2 g'_2 (\dot{w}^{k+1}_2)^{-1} g''_2 n''_2
\label{w2}
\eeqa
for some $g'_1 \in G^{k+1}_1,$ $g''_1 \in (\ol{G}^{k+1}_1)'$
$g'_2 \in (\ol{G}^{k+1}_2)',$ $g''_2 \in G^{k+1}_2,$ 
$n'_1 \in N^{k+1}_{1 \e_{k+1}},$ $n'''_1 \in \ol{N}^{k+1}_{1 +},$
$n''_1 \in \ol{N}^{k+1}_{2 -},$ and $n''_2 \in N^{k+1}_{2 \e_k}.$
The elements $w^{k+1}_i$ are minimal length representatives
in $W^k_i$ of cosets from $\ol{W}^{k+1}_i \backslash W / W^{k+1}_i.$ 
Eqs. \eqref{w1}, \eqref{w2} define them uniquely.
The elements $g^{k+1}_i \in G^{k+1}_i$ are defined by
\beqa
&& g^{k+1}_1 =
\Big( (\Th')^{-1}\Ad_{(\dot{w}_2^k \ldots \dot{w}_2^1)}^{-1} g'_2
\Big) \: g'_1,
\label{g1} \\
&& g^{k+1}_2 =
g''_2 \: \Big( \Th' \Ad_{(\dot{w}_1^k \ldots \dot{w}_1^1)}^{-1} g''_1 
         \Big)
\label{g2}
\eeqa
(see \eqref{th'} and \eqref{compute}).
The described procedure has finitely many steps $k =1, \ldots, K.$
Denote $v_i = w^K_i \ldots w^1_i,$ $i=1,2.$ For their
representatives in $N(H)$ we set
$\dot{v}_i = \dot{w}^K_i \ldots \dot{w}^1_i.$
Similarly to the proof of \leref{lemma1} one shows that
\[
(v_1, v_2) \in V
\]
and that any pair $(v_1^0, v_2^0) \in V$ can be obtained through
the described procedure. At the last step of the procedure
\[
\g^K_i = \g^{v_1, v_2}_i.      
\]

We will also need the connected subgroups of $G_i^k$
and $\ol{G}_i^k$ with Lie algebras $\g_i^k \cap (\l'_i \op \a_i)$ and
$\ol{\g}_i^k \cap \Ad_{(\dot{w}_i^{k-1} \ldots \dot{w}_i^1)}
(\l'_i \op \a_i).$
They will be denoted by $L_i^k$ and $\ol{L}_i^k.$ Let $G^{k, r}_1$ and
$G^{k, r}_2$ denote the subgroups of 
$G^{k-1}_1 \times G^{k-1}_2$ defined
by 
\beqa
G_1^{k, r} =
  &&\{ (x, y) \in L_1^k \times \ol{L}_2^k |
  \Th \pr_+(x) = \pr_-
  \Ad_{(\dot{w}_2^{k-1} \ldots \dot{w}_2^1)}^{-1}(y) \}^\ci \times
\nn \\
  && (H^\ort_1 N_{1 \e_{k-1}} \times
  \Ad_{(\dot{w}_2^{k-1} \ldots \dot{w}_2^1)}(H^\ort_2) \, \ol{N}^k_{2-} )
\label{Gkr1}
\eeqa
and
\beqa
G_2^{k, r} =
  &&\{ (x, y) \in \ol{L}_1^k \times L_2^k |
  \pr_+ \Ad_{(\dot{w}_1^{k-1} \ldots \dot{w}_1^1)}^{-1}(x)=
  \Th^{-1} \pr_-(y) \}^\ci \times
\nn \\
  && (\Ad_{(\dot{w}_2^{k-1} \ldots \dot{w}_2^1)} (H^\ort_1) \,
  \ol{N}^k_{1+} \times
  H^\ort_2 N_{1 \e_k}).
\label{Gkr2}
\eeqa
Similarly to the proof of \eqref{Gr}, using the fact $\La_i \sub H,$
one shows that the groups $G^{k, r}_i$ are connected.
They will play the role of $G^r$ (as left and right factor)
at the $k$-th step of the procedure. In particular $G^{1,r}_i = G^r,$
$i=1,2.$
\hfill \qed 
\\ \hfill \\
{\em{Proof of \thref{Grcoset}.}}
As in the proof of \thref{G+coset} we construct inductively a 
sequence of dense, open subsets of the double coset 
$[(l_1 \dot{w}_1, (\dot{w}_2)^{-1} l_2)].$ The set from the first
step consists of the following elements of $G \times G$
\beqa
&&\Big( m_+ 
        ( (\Th')^{-1} l''_+ )
        ( \pr^{-1}_+ \Th^{-1} \pr_-(l''_2) )
        ( (\Th')^{-1} l''_- )  \;
l_1 \dot{w}_1 \; l'_+ n'_+ l'_\st m'_\st l'_- n'_-, 
\nn \\
&&n''_+ l''_+ m''_\st l''_\st n''_+ l''_+ \;
(\dot{w}_2)^{-1} l_2 
        ( \Th' l'_+ )
        ( \pr^{-1}_- \Th \pr_+(l'_\st) )
        ( \Th' l'_- ) m_- \Big)
\label{hugeprod}
\eeqa
where the elements $l'_+,$ $n'_+,$ $l'_\st,$ $m'_\st,$ $l'_-,$ 
$n'_-$ belong to the factors of the following dense, open
subset of $G_+:$
\beqa
&& \left( L^1_1    \cap \Ad_{\dot{w}^1_1}^{-1}(N^1_{1+})  \right)
   \left( N^1_{1+} \cap \Ad_{\dot{w}^1_1}^{-1}(N^1_{1+})  \right)
   \left( L^1_1    \cap \Ad_{\dot{w}^1_1}^{-1}(G^1_1)     \right) \times
\nn \\
&& \times \:
   \left( M_+      \cap \Ad_{\dot{w}^1_1}^{-1}(G^1_1)     \right)
   \left( L^1_1    \cap \Ad_{\dot{w}^1_1}^{-1}(N^1_{1-})  \right)
   \left( N^1_{1+} \cap \Ad_{\dot{w}^1_1}^{-1}(N^1_{1-})  \right).
\label{product1}
\eeqa
Similarly $n''_+,$ $l''_+,$ $m''_\st,$ $l''_\st,$ $n''_+,$ 
$l''_+,$ belong to the factors of the following dense, open
subset of $G_-$
\beqa
&& \left(\ol{N}^1_{2-}  \cap \Ad_{\dot{w}^1_2}^{-1} (N^1_{2+})  \right)
   \left(L^1_2          \cap \Ad_{\dot{w}^1_2}^{-1} (N^1_{2+})  \right)
   \left(M_-            \cap \Ad_{\dot{w}^1_2}^{-1} (G^1_2)     \right)  
\times
\nn \\
&& \times \:
   \left(L^1_2          \cap \Ad_{\dot{w}^1_2}^{-1} (G^1_2)     \right) 
   \left(\ol{N}^1_{2-}  \cap \Ad_{\dot{w}^1_2}^{-1} (N^1_{2-})  \right)
   \left(L^1_2          \cap \Ad_{\dot{w}^1_2}^{-1} (N^1_{2-})  \right).
\label{product2}
\eeqa

The two sets \eqref{product1} and \eqref{product2} are constructed
according to how $\Ad_{\dot{w}_1^1}(G_+)$ intersects the dense, open
subset $N_{1+}^1 G_1^1 N_{1-}^1$ of $G$ and similarly
how $\Ad_{\dot{w}_2^1}(G_-)$ intersects $N_{2+}^1 G_2^1 N_{2-}.$ 
The set \eqref{hugeprod} is obtained by replacing the first component
of the right factor $G^r$ with \eqref{product1} and
the second component of the left factor $G^r$ with the set
\eqref{product2}.

Similarly to \eqref{2coset} one shows that the set of elements
\eqref{hugeprod}
can be represented as
\beqa
&& \Big(e, (\dot{w}^1_2)^{-1} \Big) \,
\Big( N^1_{1+}
\times(\Ad_{\dot{w}^1_2}(\ol{N}^1_{2-})\cap N^1_{2+})\Big)
\, \Big( ((\Th')^{-1} \Ad_{\dot{w}^1_2}^{-1}, \id)
(\Ad_{\dot{w}_2^1}(\ol{L}^2_1) \cap N^1_{2+}) \Big)
\times
\nn \\
&& \times \:
G_1^{2, r} \: \Big( g_1^1, g_2^1 \Big) \: G_2^{2, r} \times
\nn \\
&& \times \:
\Big( (\id, \Th'  \Ad_{\dot{w}^1_1}^{-1})
(\Ad_{\dot{w}_1^1} \ol{L}^1_1 \cap N^1_{1-}) \Big) \,
\Big( (\Ad_{\dot{w}^1_1} (\ol{N}^1_{1+})\cap N^1_{1-})
\times N_{2-}^1\Big)
\: \Big( \dot{w}^1_1, e \Big)
\label{2coset12}
\eeqa
(recall eqs. \eqref{Gkr1}, \eqref{Gkr2} for the definitions of
$G_i^{2,r}).$ One proceeds with the double -- left
$G_1^{2,r}$ and right $G_1^{2, r}$ coset in $G^1_1 \times G^1_2$
in the same way as in the proof of \thref{G+coset}.
Unfortunately the unipotent radicals of the projections
of $G_1^{k, r}$ and $G_2^{k, r}$ on the first component (and also
on the second) correspond to roots of different sign
for even $k$ and to roots of the same sign for odd $k.$
(Compare $N_{1 \e_{k+1}}^k,$ $\ol{N}_{1+}^k$ and
$\ol{N}_{2-}^k,$ $N_{2 \e_k}^k.)$ This does not allow obtaining
a closed formula for the dimensions of the cosets.
The necessary modifications will be made in the next Subsection.

At the last step of the procedure we get the double coset 
\beq
G_1^{K+1, r} (g^K_1, g^K_2) G_2^{K+1, r}
\label{laststep}
\eeq
in $G^{v_1, v_2}_1 \times G^{v_1, v_2}_2.$ The groups
$N_{i \pm}^{K+1}$ and $\ol{N}_{i \pm}^{K+1}$
are trivial and
\beqa
&&G_1^{K+1,r} = \Big( H^\ort_1 \times \Ad_{\dot{v}_2}(H^\ort_2)\Big)
\, \Big( (\id, \Ad_{\dot{v}_2} \Th')(G^{v_1, v_2}_1)' \Big) \,
Z_1^{v_1, v_2},
\label{GK+1r1} \\
&&G_2^{K+1, r} =\Big( \Ad_{\dot{v}_1}(H^\ort_1) \times H^\ort_2 \Big)
\, \Big( (\Ad_{\dot{v}_1} (\Th')^{-1}, \id)(G^{v_1, v_2}_2)' \Big) \,
Z_2^{v_1, v_2}.
\label{GK+1r2}
\eeqa
This is a consequence from  
\[
\pr_\pm( Z(L^{v_1, v_2}_i)^\ci ) = Z( L^{v_1, v_2}_i / \La_i )^\ci
\]
which can be proved as follows. The group on the left is a
connected subgroup of the one on the right and the two have one and the
same tangent Lie algebras, so they coincide.

Decompose $g^K_i \in G^{v_1, v_2}_i$ as
\[
g^K_i = f^K_i c^K_i, \;
f^K_i \in (G^{v_1, v_2}_i)', \;
c^K_i \in Z(G^{v_1, v_2}_i)^\ci.
\]
Each pair $(f^K_i, c^K_i)$ is only defined modulo
multiplication of $f^K_i$ and division of $c^K_i$ by an element from the
finite abelian group 
$(G^{v_1, v_2}_i)' \cap Z(G^{v_1, v_2}_i)^\ci,$ $i=1, 2.$
Using eqs. \eqref{GK+1r1} and \eqref{GK+1r2} one rewrites
\eqref{laststep} as
\beq
\Big( TC^{\dot{v}_1, \dot{v}_2}_{(G^{v_1, v_2})'}
(f_1^K \Ad_{\dot{v}_1} (\Th')^{-1}(f_2^K)^{-1}), e \Big) \,
\Big( (c_1^K, c_2^K) Z^{v_1, v_2} \Big) \, 
\Big( (\id, \Th' \Ad_{\dot{v}_1}^{-1})(G^{v_1, v_2}_1)' \Big)
\eeq
(recall the definition \eqref{Cv} of the group $Z^{v_1, v_2}).$

The final result of the procedure is 
\beq
[(l_1 \dot{w}_1, \dot{w}_2^{-1} l_2)] =
[(f^K_1 (\Ad_{\dot{v}_1} (\Th')^{-1} (f_2^K)^{-1}) c^K_1 \dot{v}_1,
\dot{v}_2^{-1} c^K_2)]
\label{4coset12}
\eeq
and that this coset has the following dense, open subset
\beqa
&&\Big( e, \dot{v}_2^{-1}  \Big) \, \Big( \prod_{k=1}^K E^k \Big) \times 
\nn \\
&& \times \:
\Big( TC^{\dot{v}_1, \dot{v}_2}_{(G^{v_1, v_2}_1)'}
(f_1^K \Ad_{\dot{v}_1} (\Th')^{-1}(f_2^K)^{-1}), e \Big) \, 
\Big( (c_1^K, c_2^K) Z^{v_1, v_2} \Big) 
\, \Big( (\id, \Th' \Ad_{\dot{v}_1}^{-1})(G^{v_1, v_2}_1)' \Big) \times
\nn \\
&& \times \:
\Big( \prod_{k=K}^1 F^k \Big)\, \Big( \dot{v}_1 , e \Big)
\label{3coset12}
\eeqa
where the subsets $E^k$ and $F^k$ of $G\times G$ are defined by
\beqa
E^k &=& \Big( N_{1 \e_{k+1}}^k \times
\Ad_{(\dot{w}_2^K \ldots \dot{w}_2^{k+1})}
(\Ad_{\dot{w}_2^k}(N_{2-}^k) \cap N_{2 \e_k}) \Big) \times 
\nn \\
&& \times \: \Big( ( (\Th')^{-1} \Ad_{\dot{w}_2^k}^{-1}, 
\Ad_{(\dot{w}_2^K \ldots \dot{w}_2^{k+1})})
( \Ad_{\dot{w}_2^k}(\ol{L}_2^k) \cap N_{2 \e_{k+1}}^k ) \Big),
\nn 
\eeqa
and 
\beqa
F^k &=&
\Big( (\Ad_{(\dot{w}_1^K \ldots \dot{w}_1^{k+1})},
\Th' \Ad_{\dot{w}_1^k}^{-1})
(\Ad_{\dot{w}_1^k}(\ol{L}_1^k) \cap (N_{1 \e_k}^k) ) \Big)
\times
\nn \\
&& \times \:
\Big( \Ad_{(\dot{w}_1^K \ldots \dot{w}_1^{k+1})}
(\Ad_{\dot{w}_1^k}(\ol{N}_{1+}
\cap N_{1 \e_k}^k) \times N_{2, \e_k}^k) \Big).
\nn
\eeqa
(Here as before $\e_k = (-1)^k.$)

To the coset \eqref{4coset12} one associates the pair
$(v_1, v_2) \in V,$ the $TC^{\dot{v}_1, \dot{v}_2}_{(G^{v_1, v_2})'}$
orbit of \\
$f_1^K v_1 (\Th')^{-1}(f_2^K) \in G^{v_1, v_2},$ and the
coset
\[
(c^K_1, c^K_2) Z^{v_1, v_2} \in
(Z(G^{v_1, v_2}_1)^\ci \times Z(G^{v_1, v_2}_2)^\ci)
/ Z^{v_1, v_2}.
\]
This is only defined modulo the action of the finite
group \eqref{factab} due to the freedom in the definition 
of $f^K_i$ and $c^K_i.$ 

Analogously to the proof of \thref{G+coset} one shows that two
sets of the type \eqref{3coset12} intersect if and only if they
correspond to one and the same pair $(v_1, v_2) \in V,$
orbit of $TC^{\dot{v}_1, \dot{v}_2},$
and cosets of the abelian factor group \eqref{factab}, modulo the
action of the group \eqref{cap}. The equivalence relation respects the
nonuniqueness in the definition of $c^K_i$ and $f^K_i.$ 
\hfill \qed

   {}  From this proof one also derives:
\bco{Grcos-full} If the Poisson structure on the group $G$
comes from an $r$-matrix for which $\g_\pm \sup \h$ and
the groups $\La_i$ are trivial then the $G^r$ double cosets
of $G \times G$ are classified by a pair of minimal representatives
$(v_1, v_2)$ in $W$ of cosets from $W/W_i$ and an orbit of the twisted
conjugation action of $G^{v_1, v_2}_1$ on itself defined by
\beq
\wt{TC}^{\dot{v}_1, \dot{v}_2}_g (f) =
g^{-1} f (\Ad_{\dot{v}_1} \Th^{-1} \Ad_{\dot{v}_2} \Th (g)); \:
f, g \in G^{v_1, v_2}_1.
\label{wttcr}
\eeq
{\em{(}}Simply connectedness for $G'$ is not assumed{\em{.)}}
\eco
%%%%
\subsubsection{Another view on $G^r$ double cosets and
a dimension formula}
The goal of this Subsection is to prove a formula for the dimensions
of the $G^r$ double cosets of $G \times G.$ The idea is to conjugate
the group $G^r$ by some elements of $N(H) \times N(H)$ in such a way that
the iterative procedure of the proof of \thref{Grcoset} is still
applicable to the cosets of the resulting subgroups of $G \times G.$
In the same time the corresponding groups $G^{k, r}_1$ 
and $G^{k, r}_2$ from the iteration will have unipotent radicals of their
first and
second components, of opposite sign (independent of the parity of $k).$

Let $w^\max$ and $w^\max_i$ be the elements of maximal
length in the Weyl groups $W$ and $W_i,$ $i=1,2.$
Denote with $u^\max_i$ the representatives of minimal
length in $W$ of the cosets $w^\max  W_i \in W/ W_i.$
More explicitly:
\beq
u^\max_i = w^\max (w^\max_i)^{-1}.
\label{uexpl}
\eeq
Denote 
\beqa
&&\wh{G}^r_1 = (\id \times \Ad_{\dot{u}_2^\max}) G^r
\sub G \times G,
\nn \\
&&\wh{G}^r_2 = (\Ad_{\dot{u}_1^\max} \times \id ) G^r
\sub G \times G.
\nn
\eeqa
The double -- left $\wh{G}^r_1$ and right $\wh{G}^r_2$ cosets of
$G \times G$ will be called for simplicity 
$\wh{G}^r_1-\wh{G}^r_2$ cosets  of $G \times G.$
There is a simple connection between them and the $G^r$ double cosets
of $G \times G:$ 
\beq
[( g_1, g_2)] = \Big( e, (\dot{u}^\max_2)^{-1} \Big)
\Big( \wh{G}^r_1 (g_1 (\dot{u}^\max_1)^{-1}, \dot{u}^\max_2 g_2)
\wh{G}^r_2 \Big)
\Big( \dot{u}^\max_1, e \Big).
\eeq
    
\bth{Grdim} 1. Any $\wh{G}^r_1-\wh{G}^r_2$ coset of $G \times G$
is of the form
\beq
\wh{G}^r_1 (f c_1 \dot{v}_1 (\dot{u}^\max_1)^{-1},
(\dot{v}_2 (\dot{u}^\max_2)^{-1})^{-1} c_2)
\wh{G}^r_2
\label{Gr12}
\eeq
where $v_i$ are minimal length representatives in $W$ of cosets from
$W/W_i,$ $f \in (G^{v_1, v_2}_1)',$ and
$c_i \in Z(G^{v_1, v_2}_i)^\ci.$ All such cosets are classified
by the same data as in \thref{Grcoset} and are related to the $G^r$
double cosets of $G\times G$ via
\beq
[(f c_1 \dot{v}_1, \dot{v}_2^{-1} c_2)] =
\Big( e, (\dot{u}^\max_2)^{-1} \Big) \,
\Big( \wh{G}^r_1 (f c_1 \dot{v}_1 (\dot{u}^\max_1)^{-1},
\dot{u}^\max_2 \dot{v}_2^{-1} c_2)
\wh{G}^r_2 \Big)
\, \Big( \dot{u}^\max_1, e \Big).
\label{Grcoss}
\eeq

2. The dimension of the cosets \eqref{Gr12} and \eqref{Grcoss} is 
\beq
2 \dim \g - 2 \dim \n_+ - \dim (\g^{v_1, v_2}_1)' +
l(v_1) + l(v_2) +
\dim TC^{\dot{v}_1, \dot{v}_2}_{(G^{v_1, v_2}_1)'} (f) -
\dim Z^{v_1, v_2}.
\label{G12dim}
\eeq
\eth
\proof The images of $\l_i,$ $L_i,$ $A_i,$ $\m_\pm,$ and $M_\pm$
under conjugation with $\dot{u}^\max_i$ will be denoted by
$\wh{\l}_i,$ $\wh{L}_i,$ $\wh{A}_i,$ $\wh{\m}_\mp,$ and $\wh{M}_\mp.$

Since $w^\max,$ $w^\max$ are the elements of maximal length in
$W,$ $W_i,$ they map (simple) positive roots of $\g,$ $\l_i$ into
(simple) negative roots of the same algebras.
Formula \eqref{uexpl} for $u^\max_i$ implies that
$\Ad_{\dot{u}^\max_i}$ maps

(1) (simple) positive roots of $\l_i$ into (simple) positive roots of
$\wh{\l}_i,$ and

(2) positive roots of $\g$ that are not roots of $\l_i$ into
negative roots of $\g$ that are not roots of $\wh{\l}_i.$ \\ From (1) we
get that $\wh{\l}_i$ are reductive subalgebras of $\g$
generated by $\h$ and root spaces of simple roots of $\g.$  
Property (2) implies that the Lie algebras $\wh{\m}_\pm$
are spanned by negative root spaces of $\g$ and
$(u^\max_i)^{-1}(\h^\ort_i).$

The groups $G^r_i$ can be more explicitly written as
\beqa
&&\wh{G}^r_1 = (\wh{M}_-, M_+) \,
 \{ (x, y) \in (\wh{L}'_1 \wh{A}_1) \times (L'_2 A_2) | \, 
   \Th \pr_+ \Ad_{\dot{u}^\max_1}^{-1} (x) = y \},
\nn \\
&&\wh{G}^r_2 = (M_-, \wh{M}_+)
\{ (x, y) \in (L'_1 A_1) \times (\wh{L}'_2 \wh{A}'_2) | \,
   \Th \pr_+ (x) = \pr_- \Ad_{\dot{u}^\max_2}^{-1} (y) \}.
\nn
\eeqa
The two groups are of the same type as $G^r.$ On Lie algebra level
the role of the Cayley transformation
$\th: \l'_1 \op \a_1 \ra \l'_2 \op \a_2$
is played by the map
\[
\th \Ad_{\dot{u}^\max_1}^{-1}: \wh{\l}'_1 \op \wh{\a}_1 \ra
\l'_2 \op \a_2
\]
for the first group, and by
\[
\Ad_{\dot{u}^\max_2} \th: \l'_1 \op \a_1 \ra \wh{\l}'_2 \op \wh{\a}_2
\]
for the second one.

The procedure from the proof of \thref{Grcoset} can be applied
to describe the set of $\wh{G}^r_1 - \wh{G}^r_2$ cosets  of
$G \times G.$ It starts with $\g^1_1 = \wh{\l_1},$
$\ol{\g}^1_1 = \l_1,$
$\g^1_2 = \l_2,$ and $\ol{\g}^1_2 = \wh{\l_2}.$ The final result
depends on a choice of a pair $(\wh{v}_1, \wh{v}_2)$ of minimal
length
representatives in $W$ of cosets from $W/ \wh{W}_i.$ Denote the set
of such pairs by $\wh{V}$ and recall the definition of the set $V$ from
Sect.~5.1.1. Properties (1)--(2) of the elements $u^\max_i \in W$
imply 
\[
(v_1, v_2) \in V \; \mbox{if and only if} \;
(v_1 (u^\max_1)^{-1}, v_2 (u^\max_2)^{-1}) \in \wh{V}.
\]
For any $(v_1, v_2) \in V$ we will use the representatives
$\dot{v}_1 (\dot{u}^\max_1)^{-1},$
$\dot{v}_2 (\dot{u}^\max_2)^{-1}$
of $v_1 (u^\max_1)^{-1},$ $v_2 (u^\max_2)^{-1}$ in $N(H).$
As in \thref{Grcoset} the procedure for the
$\wh{G}^r_1 - \wh{G}^r_2$ cosets of $G \times G$
gives that any such coset is equal to a
coset of the form \eqref{Gr12} for some pair
$(v_1 u^\max_1, v_2 u^\max_2) \in \wh{V}$
(i.e. $(v_1, v_2) \in V).$ 

The important point now is that at each step
of the procedure
the unipotent radicals of
the first and second components of the (new) groups $G^{k, r}_1$ and 
$G^{k, r}_2$ have opposite signs. This gives a closed formula for the
dimension of the coset \eqref{Gr12}:
\[
2 \dim \g - \dim (\g^{v_1, v_2}_1)' - l(v_1 (u^\max_1)^{-1}) -
l(v_2 (u^\max_2)^{-1}) +
\dim TC^{\dot{v}_1, \dot{v}_2}_{(G^{v_1, v_2}_1)'} (f)
- \dim Z^{v_1, v_2}.
\]
(The proof is identical to the one of the dimension formula
\eqref{dim}.)
Formula \eqref{G12dim} follows from it and
\[
l(v_i u^\max_i) = \dim \n_+ - l(v_i).                          
\]
The last equation is a consequence from properties (1)--(2) of $u^\max_i$
and the fact that the length of $w \in W$ is equal to the number
of positive roots $\al$ of $\g$ for which $w(\al)$ is a negative
root.
\hfill \qed
\subsection{The group $d(G) \cap G^r$}
Let $\Sig$ denote the group $d(G) \cap G^r.$
It is a discrete subgroup of $G \times G$ since
$d(\g) \cap \g^r = 0.$

If $(x, y) \in \Sig$ then $y=x$ and $\pr_+(x) = \pr_-(x).$
So $(x, \pr_+(x)) \in \Sig_-$ and one can define
a map
\beq
\phi: \Sig \ra \Sig_-, \; \mbox{by} \;
\phi(x, y) = (x, \pr_+(x))
\label{phi}
\eeq
(recall the notation from Sect.~4.2).
It is an isomorphism because of \eqref{Gr}.

\bpr{Sigr} For a Poisson structure on $G$ 
for which $G^r$ is a closed subgroup of $G \times G$ the groups
$\Sig$ and $\Sig_-$ are canonically isomorphic via \eqref{phi}.
For any factorizable Poisson structure on $G$ 
\beq
\Sig =i(H) \cap G^r.
\label{sigg}
\eeq
\epr

In the general case \eqref{sigg} is proved in the same way as
\thref{Sig+}. 

Note that all points of $\Sig$ are fixed under the dressing action
of $G^r$ on $G$ (i.e. they are 0 dimensional leaves of $G).$
Thus \eqref{sigg} can be viewed as a dressing analog of the 
classical fact that the 
center of a reductive group $R$ belongs to any maximal torus
of $R.$ The theorem deals with a fixed torus because the
Poisson structure depends on the choice of the torus.
%%%%%%%%%%%%%%%%%%
\subsection{Symplectic leaves of G}
As a direct consequence from the results of the previous Subsections
and Theorems \ref{t1.STS} and \ref{tint}, 
we obtain the following classification result for the
set of symplectic leaves on the Poisson--Lie group $G.$

\bth{Grsymp} Let $G$ be a complex reductive Poisson--Lie group
with simply connected derived subgroup $G'.$

1. The symplectic leaves in $d(G)/ \Sig$ are classified
by the data of \thref{Grcoset}.
The dimension of the leaf
corresponding to the pair of minimal length representatives $(v_1, v_2)$
of cosets from $W/ W_i$ and the $TC^{\dot{v}_1, \dot{v}_2}$ orbit through
$f \in (G^{v_1, v_2}_1)'$ is
\[
\dim(\l'_1 \op \a_1) + 2 \dim \h^\ort_1 - \dim (\g^{v_1, v_2}_1)' +
l(v_1) + l(v_2) +
\dim TC^{\dot{v}_1, \dot{v}_2}_{(G^{v_1, v_2}_1)'} (f_1) -
\dim Z^{v_1, v_2}.
\]

2. The symplectic leaves of $G$ are coverings of the ones
of $d(G)/\Sig$ under the projection $G \ra d(G)/\Sig.$
The group $\Sig$ can be computed explicitly from \prref{Sigr}.
\eth

In a particular (but important) case this Theorem simplifies
a lot (see Corollaries \ref{cGrcos-full} and \ref{cSig}).

\bco{Grsympfull} Let $G$ be a factorizable complex reductive 
Poisson--Lie group for which the associated $r$-matrix 
satisfies $\g_\pm \sup \h$ and the groups $\La_i$ are trivial.

Then the symplectic leaves of $d(G) / \Sig$ are classified by
a pair of minimal length representatives $(v_1, v_2)$ of cosets from $W/
W_i$ and an orbit of the twisted conjugation action
$\wt{TC}^{\dot{v}_1, \dot{v}_2}$ of $G^{v_1, v_2}_1$ on itself
defined in \eqref{wttcr}.

The dimension of the leaf that corresponds to this data is
\[
\dim \l_1 - \dim \g^{v_1, v_2}_1 +
\dim TC^{\dot{v}_1, \dot{v}_2}_{G^{v_1, v_2}_1} (f)
+ l(v_1) + l(v_2).
\]

The group $\Sig$ is explicitly given by
\[
\Sig = \{ (h, h) | \, h \in H, \, \Th(h) = h \}
\cong \ker / (1 - \th)(\ker)
\]
where $\ker$ denotes the kernel of the exponential map for the maximal
torus $H.$
\eco

\bre{comparison} Consider the special choice of the data
of \thref{Grsymp} and \coref{Grsympfull} when $v_2$ is taken the identity
element $\id$ of $W$ and only a subgroup $F_1$ of the abelian factor
group $(Z(G^{v_1, \id}_1)^\ci \times Z(G^{v_1, \id}_2)^\ci)/Z^{v_1, \id}$
is considered.
It is defined as the image of the projection of the
subtorus $Z(G^{v_1, \id}_1)^\ci \times e$ of
$Z(G^{v_1, \id}_1)^\ci \times Z(G^{v_1, \id}_2)^\ci.$ 
The simplectic leaves of $G$ corresponding
to this special choice of the data are the symplectic leaves of the
Poisson subgroup $G_-$ of $G.$ 

The connection with \thref{SymlG-} and \coref{symlG+full} is as follows.
For any minimal representative $v_1$ in $W$ of a
coset from $W/W_1$ one has $G^{v_1, \id}_1 = G^{v_1}$ and the twisted
conjugation actions $TC^{\dot{v}_1, e}$ and $TC^{\dot{v}_1}$ are 
the same. Using
\[
Z^{v_1, \id}_1 =
\{ (x, y) \in Z(L^{v_1, \id}_1)^\ci \times Z(L^{v_1, \id}_2)^\ci |
\Th \pr_+ (x) = \pr_-(y) \}^\ci
\]
one identifies the part of the data of Theorems \ref{tG+coset} and
\ref{tGrcoset} coming from the abelian factor groups
\[
F_1 \cong Z(G^{v_1})^\ci / (\Cong^{v_1}(L^{v_1}) H^\ort_1 
\Ad_{\dot{v}_1}(H^\ort_1)).
\]
The fact that the groups $L'_i$ are simply connected implies 
$\La_i \in Z(L_i)^\ci.$ From it one also gets 
\[
\La^{v_1}_1 =
d( (G^{v_1})'  \cap Z(G^{v_1})^\ci) (e, \La_1)
\sub (G^{v_1})' \times Z(G^{v_1})^\ci
\]
where $d(.)$ stays for the diagonal embedding. This gives the 
connection between the ``modulo'' parts of 
Theorems \ref{tG+coset} and \ref{tGrcoset}.

The iterative procedures from the proofs
of these Theorems are also related.
Let us start with the $G^r$ double coset of 
the element $(l_1 \dot{w}_1, e) \in G \times G$ 
and the $i_+(G_+)$ double coset of the same element
but considered in $G \times ((L'_1 \times A_1)/\La_1).$
The reductive Lie algebras $\g^k,$ $\ol{\g}^k$
and $\g_1^k,$ $\ol{\g}_1^k$ from these proofs are related by
\beqa
&&\g_1^{2k} = \g^k, 
\nn \\
&&\ol{\g}_1^{2 k-1} = \ol{\g}^k.
\nn
\eeqa
At odd (even) steps $k$ the algebras $\g_1^k$ 
(respectively $\ol{\g}_1^k)$ are kept unchanged. 
\ere

\bex{stand-Gl} If the complex semisimple group $G$ is equipped with
the standard Poisson structure then any double $G^r$
coset of $G \times G$ is equal to a coset of the form
$[(h_1 \dot{v}_1, \dot{v}_2^{-1} h_2)]$ for some $h_i \in H$
and $v_i \in W.$
The procedure from the proof of \thref{Grcoset} has only one
step and it produces the following representation
of the above coset (see also \coref{Grcos-full})
\beqa
&&\Big( \dot{v}_2^{-1}, e \Big) 
\Big( N_+ \times (\Ad_{\dot{v}_2}(N_-) \cap N_+) \Big) 
\Big( \wt{TC}^{\dot{v}_1, \dot{v}_2}_H
( h_1 \Ad_{\dot{v}_1}(h_2)), e \Big) \times 
\nn \\
&&\times \,
\Big\{(h^{-1},\Ad_{\dot{v}_1}(h)) \Big| h \in H \Big\}
\Big((\Ad_{\dot{v}_1}(N_+) \cap N_-) \times  N_- \Big)
\Big( \dot{v}_1, e \Big),
\label{Hprod}
\eeqa
where
\beq
\wt{TC}^{\dot{v}_1, \dot{v}_2}_h (h' ) = h^{-1} h'
\Ad_{\dot{v}_1 \dot{v}_2}(h).
\label{wth}
\eeq
All such double cosets and thus all symplectic leaves of
$d(G)/ \Sig$ are classified by a pair $(v_1, v_2) \in W$ and
an orbit of the action \eqref{wth}. Both \coref{Grsympfull}
and \eqref{Hprod} imply that the dimension of the corresponding
symplectic leaf is
\[
\dim \wt{TC}^{\dot{v}_1, \dot{v}_2}_H(h) + l(v_1) + l(v_2)
\]
The group $\Sig$ is isomorphic to the group $\Sig_-$ from \exref{Sigstand}
(see \eqref{SigmaSt1} and \eqref{SigmaSt2}).
This Example agrees with the results in \cite{HL, HKKR}.
\eex
\bex{CG-Gl}Here we discuss the set of symplectic leaves
of the group $SL(n+1)$ equipped with Cremmer--Gervais
Poisson structure.

The first step is the classification of the $SL(n+1)^r$
double cosets of $\D(SL(n+1)).$ Recall from \exref{CGcosets} that the
minimal
length representatives in $W \cong S^{n+1}$ of the cosets from
$W/W_1 \cong  S^{n+1}/S_1^n$ are the permutations
$\sig^j_1,$ $j=0, \ldots, n$ defined in \eqref{sig1}.
In the same way one finds that the
minimal length representatives for the cosets from 
$W/W_2  \cong S^{n+1}/S^n_2$ are the permutations
\[
\sig^k_2 = (1 \ldots (n+1-k))^{-1}, \; k=0, \ldots, n.
\]
As in \exref{CGcosets} we will assume that all elements of the
Weyl group of $SL(n+1)$ are represented in $N(H)$
by multiples of permutation matrices.

(1) Case $j+k \geq n.$  The stable subgroup 
$G^{\dot{\sig}^j_1, \dot{\sig}^k_2}_1$ is given by
\[
G^{\dot{\sig}^j_1, \dot{\sig}^k_2}_1 =
\{ \diag(a_1, \ldots, a_{n-k}, A, a_{j+1}, \ldots, a_{n+1}) \, |
\, a_p \in \Cset, \, A \in GL(j+k -n), \det A \prod_p a_p =1 \}.
\]
Under the twisting conjugation action 
$TC^{\dot{\sig}^j_1, \dot{\sig}^k_2}$ of this this group
on itself the element
$\diag(a_1, \ldots, a_{n-k}, A, a_{j+1}, \ldots, a_{n+1})$
acts as follows
\beqa
&&\diag(b_1, \ldots, b_{n-k}, B, b_{j+1}, \ldots, b_{n+1}) \mapsto
\label{bBb} \\
&&\diag(a_1^{-1} b_1 a_{\mu(1)}, \ldots, 
a_{n-k}^{-1} b_{n-k} a_{\mu(n-k)}, A^{-1} B A, 
a_{j+1}^{-1} b_{j+1} a_{\mu(j+1)}, \ldots, 
a_{n+1}^{-1} b_{n+1} a_{\mu(n+1)})
\nn
\eeqa
where the permutation $\mu \in S^{n+1}$ is given by
\[
\mu = ( (j+1) \ldots (n+1) \, (n-k) \ldots 1).
\]
There is a one to one correspondence between the set of 
$TC^{\dot{\sig}^j_1, \dot{\sig}^k_2}$ orbits on
$G^{\dot{\sig}^j_1, \dot{\sig}^k_2}_1$ and the conjugation
orbits of $GL(j+k-n).$ To the orbit of the element
\eqref{bBb} of $G^{\dot{\sig}^j_1, \dot{\sig}^k_2}_1$
is associated the conjugation orbit of $B \in GL(j+k-n).$

Therefore the $SL(n+1)^r$ double cosets of 
$SL(n+1)^{\times 2} \cong \D(SL(n+1))$
and the symplectic leaves of $SL(n+1) / \Sig$ 
corresponding to the permutations $\sig^j_1,$ $\sig^k_2$ with
$j+k \geq n$
are classified by the conjugation orbits of $GL(j+k-n)$ on itself.
The dimension of the leaf corresponding of the orbit 
trough $B \in GL(j+k-n)$ is
\[
(2n-j-k)(j+k+1) + \dim \{ \Ad_A B | \, A \in GL(j+k -n) \}.
\]
The group $\Sig$ coincides with the group $\Sig_-$ from
\exref{CGcosets} (see \eqref{SigmaCG}) and the symplectic
leaves in $SL(n+1)$ are coverings of the
ones of $SL(n+1)/\Sig$ under the covering map 
$SL(n+1) \ra SL(n+1)/\Sig.$

(2) Case $j+k \leq n-1.$ Then the stable subgroup
$G^{\dot{\sig}^j_1, \dot{\sig}^k_2}_1$ is
\[
\{ \diag(a_1, \ldots, a_{j+1}, A, a_{n-k+1}, \ldots, a_{n+1}) \, |
\, a_p \in \Cset, \, A \in GL(n-j-k-1), \det A\prod_p a_p = 1 \}.
\]
The $TC^{\dot{\sig}^j_1, \dot{\sig}^k_2}$ action 
of the element
$\diag(a_1, \ldots, a_{j+1}, A, a_{n-k+1}, \ldots, a_{n+1})$
on \\ $\diag(b_1, \ldots, b_{j+1}, B, b_{n-k+1}, \ldots, b_{n+1}) $
gives the element
\[
\diag(a_1^{-1} b_1 a_{\mu(1)}, \ldots,
a_{j+1}^{-1} b_{j+1} a_{\mu(j+1)}, A^{-1} B A,
a_{n-k+1}^{-1} b_{n-k+1} a_{\mu(n-k+1)}, \ldots,
a_{n+1}^{-1} b_{n+1} a_{\mu(n+1)})
\]
where
\[
\mu = (1 \ldots (j+1))^{-1} \, ((n-k+1) \ldots (n+1)) \in S^{n+1}.
\]
The symplectic leaves of $SL(n+1) / \Sig$
corresponding to $\sig^j_1,$ $\sig^k_2$ in the case $j+k \leq n-1$
are classified by the conjugation orbits of $GL(n-j-k-1)$ on itself
and a complex nonzero number. (It comes from the fact that
the product of the entries $a_1,$ $\ldots,$ $a_{j+1}$ is invariant under
the twisting conjugation action.)
The dimension of the leaf corresponding to the orbit
through $B \in GL(n-j-k+1)$ is
\[
(2n-j-k-2)(j+k+1) + 2 n + 
\dim \{ \Ad_A B | \, A \in GL(n-j-k-1) \}.
\]

The group $\Sig$ is the same as in case (1) and the symplectic leaves
of $SL(n+1)$ are recovered in the same way.
\eex
%%%%%%%%%%%%%%%%% References %%%%%%%%%%%%%%%%%%%%%%%%%%%%%%%%%%%%%%%%%%%%
    
%%%%%%%%%%%%%%%%%%%%%%%%%%%%%%%%%%%%%%%%%%%%%%%%%%%%%%%%%%%%%%%%%%%%%%%%%
%%%%%%%%%%%%%%%%%%%%%%%%%%%%%%%%%%%%%%%%%%%%%%%%%%%%%%%%%%%%%%%%%%%%%%%%%
\end{document}